\documentclass[11pt]{article}
\usepackage[english]{babel}
\usepackage[T1]{fontenc}
\usepackage{verbatim}
\usepackage{amsmath}
\usepackage{amssymb}
\usepackage{graphicx}

\newcommand{\R}{\mathbb R}

\newcommand{\C}{   {\ifmmode{{\mathbb C}}\else{$\mathbb C$}\fi}}

\newcommand{\norm}[1]{\left\Vert #1\right\Vert}

\newcommand{\Z}{\mathcal{Z}}

\newcommand{\cci}{\mathcal{C}^{\infty}}
\newcommand{\D}{\mathcal{D}}

\newcommand{\Li}{\mathcal{L}}
\newcommand{\Lo}{\mathcal{L}_0}

\newcommand{\n}{\mathfrak{N}}

\newcommand{\X}{\mathcal{X}}

\newcommand{\ab}{\mathfrak{a}}
\newcommand{\h}{\mathfrak{h}}
\newcommand{\g}{\mathfrak{g}}
\newcommand{\z}{\mathfrak{z}}

\newcommand{\dt}{\frac{d}{dt}}
\newcommand{\dto}{\frac{d}{dt}_{|t=0}}
\newcommand{\dso}{\frac{d}{ds}_{|s=0}}
\newcommand{\ddx}[1]{\frac{\partial}{\partial #1}}

\newcommand{\Span}{\mbox{Span}}
\newcommand{\ad}{\mbox{ad}}
\newcommand{\Ad}{\mbox{Ad}}
\newcommand{\rank}{\mbox{rank}\ }

\newcommand{\ds}{\displaystyle}
\newcommand{\demo}{\noindent\textit{Proof. }}

\begin{document}

\title{\bf Almost-Riemannian Geometry\\ on Lie Groups}

\author{Victor AYALA\footnote{Universidad Cat\'olica del Norte, Antofagasta Chile, E-mail: vayala@ucn.cl, Partially supported by Proyecto Fondecyt No 1150292 and Pro Fondecyt VRIDT 2014 U. Cat\'olica del Norte.}, Philippe JOUAN\footnote{Lab. R.~Salem, CNRS UMR 6085, Universit\'e
    de Rouen, avenue de l'universit\'e BP 12, 76801
    Saint-\'Etienne-du-Rouvray France. E-mail: Philippe.Jouan@univ-rouen.fr}}

\date{\today}

\maketitle
\begin{abstract}
A simple Almost-Riemannian Structure on a Lie group $G$ is defined by a linear vector field (that is an infinitesimal automorphism) and $\dim(G)-1$ left-invariant ones. We state results about the singular locus, the abnormal extremals and the desingularization of such ARS's, and these results are illustrated by examples on the 2D affine and the Heisenberg groups.

These ARS's are extended in two ways to homogeneous spaces, and a necessary and sufficient condition for an ARS on a manifold to be equivalent to a general ARS on a homogeneous space is stated.

\vskip 0.2cm

Keywords: Lie groups; Linear vector fields; Almost-Riemannian geometry.

\vskip 0.2cm
\end{abstract} 







\section{Introduction}

An almost-Riemannian structure (ARS in short) on an $n$-dimensional differential manifold can be defined, at least locally, by a set of $n$ vector fields, considered as an orthonormal frame, that degenerate on some singular set
(it can also be defined by a non regular full rank  distribution and a metric, see Section \ref{BD}). 
This geometry goes back to \cite{Grushin} and \cite{Takasu}. It appears as a part of sub-Riemannian geometry, and has aroused some interest, as shown by the recent papers \cite{ABS08}, \cite{ABCGS10}, \cite{BCST09}, \cite{BCGJ11}, \cite{BS08}, \cite{BCGM14}, \cite{BCGS}, \cite{BCG13}, \cite{BL}.

On the other hand the (improperly) so-called linear systems on Lie groups have been studied quite thoroughly for ten years, in particular by the authors (see \cite{AH97}, \cite{AK03}, \cite{ARSM09}, \cite{AS01}, \cite{AT99} and  \cite{Jou08}, \cite{Jou09}, \cite{Jou11}, \cite{Jou12}, \cite{DJ14}).

It is very natural to define an ARS on an $n$-dimensional Lie group by $n$ left-invariant or affine vector fields, the rank of which is equal to $n$ on a proper open and dense subset and that satisfy the rank condition.

As the invariant vector fields can be projected to homogeneous spaces, and affine vector fields can be defined on homogeneous spaces, the previous definition of an ARS can be extended to homogeneous spaces. In Section \ref{EME} it is proved (under some technical assumptions) that an ARS that generates a finite dimensional Lie algebra is equivalent (globally or locally, according to the technical assumptions) to an ARS on a Lie group or an homogeneous space as previously defined (see Theorem \ref{ARSequivalence} and Corollary \ref{ARSLieequi} of Section \ref{EME}). The ARS's we define on Lie groups and homogeneous spaces thus appear as models for an interesting (at least in our opinion) class of ARS's.

Excepted in Section \ref{EME}, devoted to equivalence, we restrict ourselves to what we call \textit{simple ARS's}, that is ARS's defined by one linear vector field and $n-1$ invariant ones. Notice that the famous Grushin plane is of that type, on the Abelian Lie group $\R^2$.

After having stated some basic definitions and facts (Section \ref{BD}) we turn our attention to the singular locus, that is the set of points where the vector fields fail to be independant. It is an analytic set, but not a subgroup, not even a submanifold, in general. In Section \ref{SL} sufficient conditions for the singular locus to be a submanifold or a subgroup are stated (Theorems \ref{Zvariete} and \ref{Zgroupe}). This locus is very important in what concern the structure of ARS's, in particular in view of a classification. The study is illustrated by many examples, in particular on the Heisenberg group, which show that the ARS's on Lie groups, even the simple ones, go much farther than the 2D or 3D generic case (see \cite{ABS08}, \cite{BCGM14}).

The Hamiltonian equations of the PMP are computed in Section \ref{HE}. These computations are rather standard, but more complicated than in the invariant case, because the Hamiltonians can here depend on the state. They allow to get a complete characterization of the abnormals, it is the object of Theorem \ref{Anormales} in Section \ref{Abnormal}.

Before providing examples we show that simple ARS's can be globally desingularized in a standard way: we establish in Section \ref{ARSSRS} a relation between a simple ARS and a codimension one regular sub-Riemannian structure on the semi-direct product of the Lie group with the real line.

In Section \ref{Examples} two examples are computed, the first one on the two-dimensional solvable group, and the other one on the Heisenberg group.

The computation of some diffentials and power series is postponed to the appendix.

In order to avoid to lengthen the paper, we have chosen to not recall the basic facts about optimal control and sub-Riemannian geometry. The reader is referred to textbooks as \cite{PBGM62}, \cite{ABB14}, \cite{AS04}, \cite{Jurdjevic97}.


\section{Basic definitions}\label{BD}

\subsection{Linear vector fields}

In this section the definition of linear vector fields and some of their properties are recalled. More details can found in \cite{Jou09}.

Let $G$ be 
a connected Lie group and $\g$ its Lie algebra (the set of left-invariant vector fields, identified with the tangent space at the identity). A vector field on $G$ is said to be {\it linear if its flow is a one parameter group of automorphisms}. Notice that a linear vector field is consequently analytic and complete.

The flow of a linear vector field $\X$ will be denoted by $(\varphi_t)_{t\in\R}$.

The following characterizations will be useful in the sequel.

\vskip 0.2cm
\noindent \textit{\textbf{Characterization of linear vector fields}}

{\it Let $\X$ be a vector field on a connected Lie group $G$. The following
conditions are equivalent:
\begin{enumerate}
\item
$\X$ is linear;
\item
$\X$ belongs to the normalizer of $\g$ in the algebra $V^{\omega}(G)$ of analytic vector fields of $G$,
that is
\begin{equation} \label{deriv}
\forall Y\in \g \qquad [\X,Y]\in\g
\end{equation}
and verifies $\X(e)=0$;
\item
$\X$ verifies
\begin{equation}
 \forall g,g'\in G \qquad \X_{gg'}=TL_g.\X_{g'}+TR_{g'}.\X_g \label{Bourbak}
\end{equation}
\end{enumerate}}

According to (\ref{deriv}) one can associate to a given linear vector field $\X$ the derivation $D$ of $\g$ defined by:
$$
\forall Y\in \g \qquad DY=-[\X,Y],
$$
that is $D=-\ad(\X)$.
The minus sign in this definition comes from the formula
$[Ax,b]=-Ab$ in $\R^n$. It also enables to avoid a minus sign in the useful formula:
\begin{equation}\label{equationfonda}
\forall Y\in\g, \quad \forall t\in \R \qquad \varphi_t(\exp
Y)=\exp(e^{tD}Y).
\end{equation}

In the case where this derivation is inner, that is $D=-\ad(X)$ for some left-invariant vector field $X$ on $G$, the linear field splits into $\X=X+\mathcal{I}_*X$, where $\mathcal{I}$ stands for the diffeomorphism $g\in G\longmapsto \mathcal{I}(g)=g^{-1}$. Thus $\X$ is the sum of the left-invariant vector field $X$ and the right-invariant one $\mathcal{I}_*X$. The flow of such a vector field is given by
\begin{equation}
 \varphi_t(g)=\exp(-tX)g\exp(tX) \label{flot}
\end{equation}

An \textit{affine vector field} is a element of the normalizer $\n$ of $\g$ in $V^{\omega}(G)$, that is
$$
\n=\mbox{norm}_{V^{\omega}(G)}\g=\{F\in V^{\omega}(G);\ \forall Y\in \g, \quad
[F,Y]\in\g\},
$$
so that an affine vector field is linear if and only if it vanishes at the identity.

It can be shown (see \cite{AT99} or \cite{Jou09}) that an affine vector field can be uniquely decomposed
into a sum
$$
F=\X+Z
$$
where $\X$ is linear and $Z$ right-invariant.


\subsection{Almost-Riemannian structures}
For all that concern general sub-Riemannian geometry, including almost-Riemannian one, the reader is referred to \cite{ABB14}.

\newtheorem{DefARS}{Definition}
\begin{DefARS} \label{DefARS}
An almost-Riemannian structure on a smooth $n$-dimensional manifold $M$ is a triple $(E,f,\left\langle .,.\right\rangle)$ where:
\begin{enumerate}
	\item $E$ is a rank $n$ vector bundle on $M$;
	\item $f:\ E\longmapsto TM$ is a morphism of vector bundles;
	\item $(E,\left\langle .,.\right\rangle)$ is an Euclidean bundle, that is $\left\langle .,.\right\rangle_q$ is an inner product on the fiber $E_q$ of $E$, smoothly varying w.r.t. $q$;
\end{enumerate}
assumed to satisfy the following properties:
\begin{enumerate}
	\item[(i)] The set of points $q\in M$ such that the restriction of $f$ to $E_q$ is onto is a proper open and dense subset of $M$;
	\item[(ii)] The modulus $\overline{\varXi}$ of vector fields of $M$, defined as the image by $f$ of the modulus of smooth sections of $E$ satisfies the rank condition.
\end{enumerate}
\end{DefARS}

\noindent\textbf{Remarks}
\begin{enumerate}
	\item The set of points of $M$ where the rank of $f(E_q)=\varXi_q$ is less than $n$ is called the singular locus of the ARS and will be denoted by $\Z$ in the sequel. If $M\setminus \Z$ was not required to be proper (see (i)), that is if $\Z$ could be empty, then the structure could be Riemannian.
	\item This structure is trivializable if $(E,\left\langle .,.\right\rangle)$ is isomorphic to the trivial Euclidean bundle $M\times \R^n$. In that case we can choose an orthonormal basis $(e_1,\dots, e_n)$ on $\R^n$, and define $n$ vector fields on $M$ by $f_i(q)=f(q,e_i)$, $i=1,\dots,n$. The set $(f_1,\dots, f_n)$ is an orthonormal frame on $M\setminus\Z$.
\end{enumerate}

\noindent\textbf{Norm}

The almost-Riemannian norm on $\varXi_q$ is defined by
$$
\norm{v}=\min\{\norm{u};\ u\in E_q \ \mbox{ and }\ f(u)=v\}.
$$
Notice that thanks to the linearity of $f$ on the fibers, the minimum is well defined. In the trivializable case one has:
$$
\ds \norm{v}=\min\left\{\sqrt{\sum_1^n u_i^2};\ u_1f_1(q)+\dots+u_nf_n(q)=v\right\}.
$$

\vskip 0.3cm

\noindent\textbf{Admissible curves}

A curve $\gamma:\ [0,T]\longmapsto M$ is admissible if it is Lipschitz continuous (w.r.t the differential structure of $M$) and if there exists a measurable bounded function $t\longmapsto u(t)$ from $[0,T]$ into $E$ such that $\dot{\gamma}(t)=f(u(t))$ for a.e. $t\in[0,T]$. Its length is:
$$
l(\gamma)=\int_0^T\norm{\dot{\gamma}(t)} dt
$$
In the trivializable case, assuming that $\norm{u(t)}$ realizes the minimum of $\norm{\dot{\gamma}(t)}$, we get
$$
\ds l(\gamma)=\int_0^T\sqrt{\sum_1^n u_i^2(t)} dt.
$$
The almost-Riemannian distance is define as the infimum of the lengths of the admissible curves, and it is a well-known fact that this distance is finite, continuous, and induces the manifold topology (the rank condition is here essential).


\subsection{Simple ARS's}

In that paper we mainly deal with what we call \textit{simple ARS's}. They are ARS's defined on a connected Lie group $G$ by a set of $n$ vector fields
$$
\{\X,Y_1,\dots,Y_{n-1}\}
$$
where
\begin{enumerate}
	\item[(i)] $\X$ is linear;
	\item[(iii)] $Y_1,\dots,Y_{n-1}$ are left-invariant;
	\item[(iii)] $n=\dim G$ and the rank of $\X,Y_1,\dots,Y_{n-1}$ is full on a non empty subset of $G$;
	\item[(iv)] the set $\{\X,Y_1,\dots,Y_{n-1}\}$ satisfies the rank condition.
\end{enumerate}

The set where the rank of $\X,Y_1,\dots,Y_{n-1}$ is not full will be refered to as the singular locus and denoted by $\Z$.

\vskip 0.2cm

\noindent \textbf{Remarks}.
\begin{enumerate}
	\item The singular locus $\Z$ is an analytic subset of $G$. By Assumption (iii) it is not equal to $G$, and by analycity its interior is empty. On the other hand $\X(e)=0$ and it cannot be empty. Finally $G\setminus \Z$ is an open, dense and proper subset of $G$.
	\item The rank condition implies
	\begin{equation}
\mathcal{LA}\{\X,Y_1,\dots,Y_{n-1}\}=\R\X\oplus \g. \label{ARS33}
\end{equation}
\end{enumerate}

\noindent \textbf{Necessary conditions for the rank condition}

Let us denote by $\Delta$ the vector subspace of $\g$ generated by $Y_1,\dots,Y_{n-1}$, and by $D$, as usual, the derivation asociated to $\X$. If $[\Delta,\Delta]\subseteq\Delta$ and $D(\Delta)\subseteq\Delta$, then the Lie algebra generated by $\X,Y_1,\dots,Y_{n-1}$ is  equal to $\R\X\oplus \Delta$. But the rank of that Lie algebra is not full on $\Z$.
Consequently the rank condition implies that at least one of the following conditions hold:
\begin{enumerate}
	\item[(i)] $[\Delta,\Delta]\nsubseteq\Delta$ 
	\item[(ii)] $D(\Delta)\nsubseteq\Delta$
\end{enumerate}
Notice that in all cases, the full rank is obtained after one step.

\vskip 0.3cm

\noindent \textbf{Example. The Grushin plane}

The Grushin's metric on the plane corresponds to the frame $\X=(0,x_1)^T$ and $Y_1=(1,0)^T$, that is to the control system:
$$
\dot{x}=v\X+u_1Y_1 \quad\mbox{or, in coordinates,}\ \ 
\left\{
\begin{array}{ll}
\dot{x}_1 & =u_1\\
\dot{x}_2 & =vx_1
\end{array}
\right.
$$

The state space $\R^2$ is an Abelian Lie group, the field $\X$ is linear and the field $Y_1$ invariant, so that the Grushin plane is a simple ARS. The singular locus is here the line $\{x_1=0\}$.


\subsection{Notations}
In the sequel the following notations will be used:
\begin{enumerate}
\item Let $Y\in\g$. Then $Y_g$ will stand for $TL_g.Y$, where $L_g$ is the left translation by $g$, and $TL_g$ its differential.
\item To the linear vector field $\X$ we associate $F(g)=TL_{g^{-1}}.\X_g\in\g$. We also write $F_g$ for $F(g)$ when it is more convenient.
\end{enumerate}


\section{Singular locus}\label{SL}
\subsection{Introduction}

This section is devoted to the properties of the singular locus $\Z$ of
the simple ARS determined by the orthonormal frame $\{\X,Y_1,\dots,Y_{n-1}\}$. The metric does not matter here and instead of $Y_1,\dots,Y_{n-1}$ we will mainly consider the subspace $\Delta$ of $\g$ they generate.

In the Abelian simply connected case, that is when $G=\R^n$, the singular locus is easily seen to be a codimension one subspace of $\R^n$. It is natural to ask whether $\Z$ is a subgroup of $G$ in general, but the examples of Section \ref{SLE}, in particular on the group Heisenberg, show that it is not the case. Actually $\Z$ is not even a manifold in general, but only an analytic set.

By definition 
\begin{equation*}
\Z=\{g\in G\text{ }/\text{ }\rank(\X(g),Y_{1}(g),Y_{2}(g),\ldots ,Y_{n-1}(g))=n-1\}.
\end{equation*}

However we need to characterize $\Z$ in a more handable way. The hyperplane $\Delta $ can be defined as the kernel of a one-form $\omega$ to which we can associate a left-invariant one-form likewise denoted by $\omega$.

We recall from the end of the previous section that the mapping $F$ from $G$ to $\g$ is defined by
$$
F(g)=TL_{g^{-1}}.\X_g
$$
The power series expansion of $F$ and its differential are computed in the appendix (Propositions \ref{A1}, \ref{A2} and Corollary \ref{C1}).
Thanks to $\omega$ and $F$, the singular locus can be defined in a form that allows differential calculus:
$$
g\in \Z \Longleftrightarrow \X_g\in \Delta_g
        \Longleftrightarrow F_g\in\Delta
        \Longleftrightarrow \langle\omega, F_g\rangle=0.
$$
In the sequel we denote by $\psi$ the analytic function from $G$ to $\R$ defined by $\psi(g)=\langle\omega, F_g\rangle$. The singular locus is consequently the set of zeros of $\psi$.


\subsection{A sufficient condition for $\Z$ to be a submanifold}
In this section it is shown that $\Z$ is a submanifold (embedded and analytic) as soon as $\Delta$ is a subalgebra of $\g$.
\newtheorem{decompositionDez}{Lemma}
\begin{decompositionDez}\label{decompositionDez}
Let us assume that $\Delta$ is a subalgebra of $\g$. Then $\g$ can be decomposed as follows:
\begin{equation}\label{formuledecompo}
\g=\R Y_n\oplus W \oplus \R \xi
\end{equation}
where
\begin{enumerate}
\item $W$ is a $(n-2)$-dimensional subalgebra;
\item $\ds \Delta=W \oplus \R \xi$;
\item $\ds D^{-1}(\Delta) =\R Y_n\oplus W$;
\item $\langle\omega,Y_n\rangle=1$ and $D\xi=Y_n$ $(\mbox{mod }{\Delta})$.
\end{enumerate}
\end{decompositionDez}

\demo
As $\Delta$ is a subalgebra, the rank condition implies that $D^*\omega=\omega\circ D$ is different from $\omega$ and does not vanish (otherwise we would have $D\Delta\subset \Delta$ and the rank condition would not hold, see Section \ref{BD}). Consequently $D^{-1}\Delta=\ker(D^*\omega)$ is a $(n-1)$-dimensional subspace of $\g$ different from $\Delta$, and
$$
W=\Delta\cap D^{-1}\Delta
$$
is a $(n-2)$-dimensional subspace of $\g$ which turns out to be a subalgebra. Indeed let $X_1,X_2\in W$. Since $\Delta$ is a subalgebra we have $[X_1,X_2]\in\Delta$. By definition of $W$, $DX_1$ and $DX_2$ belong also to $\Delta$ so that
$$
D[X_1,X_2]=[DX_1,X_2]+[X_1,DX_2]\in \Delta
$$
In conclusion $[X_1,X_2]$ belongs to $D^{-1}\Delta$ hence to $W$ which is thus a subalgebra of $\g$.

Let $Y_n\in\g$ such that $\ds D^{-1}(\Delta) =\R Y_n\oplus W$. As $Y_n\notin \Delta$,  $\langle\omega,Y_n\rangle$ does not vanish and can be assumed to be equal to $1$.

To finish $D\Delta \nsubseteq \Delta$ hence $D\xi\notin \Delta$, and $D\xi=aY_n$ $(\mbox{mod }{\Delta})$ with $a\neq 0$. Up to a modification of $\xi$ we can choose $a=1$.

\hfill $\Box$

\newtheorem{Zvariete}{Theorem}
\begin{Zvariete} \label{Zvariete}
If $\Delta$ is a subalgebra of $\g$ then the singular locus $\Z$ is an analytic, embedded, codimension one submanifold of $G$.

Its tangent space at the identity is $\ker(D^*\omega)=D^{-1}\Delta$.
\end{Zvariete}
\demo

As $\psi$ is analytic, we have only to show that its differential does not vanish on $\Z$. According to Corollary \ref{C1} of the appendix it is given by:
$$
\forall g\in G,\ \forall Y\in\g \quad T_g\psi.Y_g=\langle\omega, DY+[F_g,Y]\rangle
$$
Consider $\xi$ as defined in Lemma \ref{decompositionDez}. Since  $\xi\in \Delta$, $F_g\in\Delta$ if $g\in\Z$, and $\Delta$ is a subalgebra, we have $[F_g,\xi]\in \Delta$ for all $g\in\Z$. On the other hand $D\xi\notin \Delta$ and:
$$
\forall g\in\Delta \qquad T_g\psi.\xi_g= \langle\omega,D\xi+[F_g,\xi]\rangle=1.
$$
This shows that $\psi$ is a submersion at all points of $\Z$ which is thus a codimension one submanifold of $G$. Moreover $F_e=0$ so that
$$
\forall Y\in\g \qquad T_e\psi.Y=\langle\omega, DY\rangle=\langle D^*\omega, Y\rangle.
$$
\hfill $\Box$
\vskip 0.2cm
\noindent\textbf{Remark}. The previous formula holds at every point where $F$ vanishes, that is on the set of singularities of $\X$.


\subsection{Sufficient conditions for $\Z$ to be a subgroup}

In order to state sufficient conditions for $\Z$ to be a subgroup we begin by a lemma which is not easily checkable in practice but which is essential in the proofs of the forthcoming Proposition \ref{Omegaclosed} and Theorem \ref{Zgroupe}.

\newtheorem{lemmeAd}[decompositionDez]{Lemma}
\begin{lemmeAd} \label{lemmeAd}
The singular locus $\Z$ is a subgroup of $G$ if and only if
\begin{equation}\label{FormulaAd}
\forall g,g'\in\ \Z\qquad \langle\omega, \Ad(g')F_g\rangle=0
\end{equation}
that is if and only if for all $g,g'\in\Z$, $ \Ad(g')F_g\in\Delta$.
\end{lemmeAd}

\demo
Let $g,g'\in G$. Then, according to Formula (\ref{Bourbak}) of Section \ref{BD},
$$
\begin{array}{ll}
F(g'g) & =TL_{(g'g)^{-1}}\X_{g'g}\\
 & =TL_{g^{-1}}TL_{(g')^{-1}}(TL_{g'}\X_g+TR_g\X_{g'})\\
  & =F(g)+\Ad(g^{-1})F(g').
\end{array}
$$
Let us assume that $g,g'\in \Z$, hence that $F(g)$ and $F(g')$ belong to $\Delta$. Then
$$
g'g\in \Z\Longleftrightarrow F(g'g)\in\Delta \Longleftrightarrow \Ad(g^{-1})F(g')\in\Delta.
$$
In the same way we get $0=F(gg^{-1})=F(g^{-1})+\Ad(g)F(g)$, so that $F(g^{-1})=-\Ad(g)F_g$, and $g^{-1}$ belongs to $\Z$ if and only if $\Ad(g)F_g\in\Delta$.

Finally we see that $\Z$ is a subgroup of $G$ if and only if Formula (\ref{FormulaAd}) holds.

\hfill $\Box$

\newtheorem{Omegaclosed}{Proposition}
\begin{Omegaclosed} \label{Omegaclosed}
The following conditions are equivalent: (i) $\Delta$ is an ideal of $\g$, (ii) the derived algebra $\D^1\g$ is included in $\Delta$, (iii)  the left-invariant one-form $\omega$ is closed.

When these conditions are satisfied the singular locus $\Z$ is a Lie subgroup of $G$ whose Lie algebra $\z$ is equal to $\ker(D^*\omega)=D^{-1}\Delta$.
\end{Omegaclosed}

\demo

The first thing to notice is that $\Z$ being closed is a Lie subgroup of $G$ as soon as it is a subgroup.

It is well-known that for a left-invariant one-form $\omega$ the Maurer-Cartan formula writes
$$
\forall X,Y\in\g \qquad d\omega(X,Y)=-\omega([X,Y])
$$
so that $d\omega$ vanishes at $e$ (hence everywhere since it is left-invariant) if and only if the derived subalgebra $\D^1\g$ is included in $\ker(\omega)$. If $Z\in\g$ is such that $\g=\Delta+\R Z$ then:
$$
\D^1\g=[\Delta,\Delta]+[\Delta,Z].
$$
It is clear that $\D^1\g\subset \Delta$ if $\Delta$ is an ideal of $\g$. Conversely if $\Delta$ is not an ideal then either $\Delta$ is not a subalgebra or $[\Delta,Z]$ is not included in $\Delta$. In both cases $\D^1\g$ is not included in $\Delta$.

\vskip 0.2cm

Let us assume $\Delta$ to be an ideal of $\g$, and let $g,g'\in\Z$. The group $G$ being connected, there exist $X_1,\dots,X_k\in\g$ such that $g=\exp(X_1)\dots\exp(X_k)$, and:
$$
\begin{array}{ll}
\Ad(g) & =\Ad(\exp(X_1)\dots\exp(X_k))\\
      & =\Ad(\exp(X_1))\dots \Ad(\exp(X_k))\\
      & =e^{\ad(X_1)}\dots e^{\ad(X_k)}.
\end{array}
$$
Each of the $e^{\ad(X_i)}$'s sends $\Delta$ into $\Delta$. Indeed for any $Y\in\g$
$$
e^{\ad(X_i)}Y=Y+\sum_{j=1}^{+\infty}\frac{1}{j!}\ad^j(X_i)Y=Y\ \  (\mbox{mod }{\D^1\g})
$$
and $e^{\ad(X_i)}Y\in\Delta$ if and only if $Y\in\Delta$.

This shows that $\Ad(g)$ sends $\Delta$ into $\Delta$, and since $F(g')\in\Delta$ we obtain $\Ad(g)F(g')\in\Delta$. According to Lemma \ref{lemmeAd} the singular locus $\Z$ is a subgroup of $G$.

The Lie algebra of $\Z$ is its tangent space at $e$. Since $\Delta$ is an ideal we know by Theorem \ref{Zvariete} that $\z=T_e\Z=\ker(D^*\omega)$.

\hfill $\Box$

\vskip 0.2cm
\newtheorem{Zgroupe}[Zvariete]{Theorem}
\begin{Zgroupe} \label{Zgroupe}
If the Lie algebra $\g$ is solvable, and $\Delta$ is a subalgebra of $\g$, then the singular locus $\Z$ is a codimension one subgroup of $G$ whose Lie algebra is $\z=D^{-1}\Delta$.
\end{Zgroupe}

\demo
\vskip 0.2cm
Notice first that according to Theorem \ref{Zvariete} the singular locus $\Z$ is a codimension one submanifold of $G$ whose tangent space at the identity is $\z=D^{-1}\Delta$. Therefore we have only to prove that $\Z$ is a subgroup.
\begin{enumerate}
\item To begin with let $\delta$ be a codimension one subalgebra of a nilpotent Lie algebra $\h$, and let us show that $\delta$ is an ideal of $\h$. Let $Z\in \h$ such that $\h=\delta+\R Z$; then its derived algebra is $
\D^1\h=[\delta,\delta]+[\delta,Z]
$.
If $\delta$ is a subalgebra but not an ideal of $\h$ then $[\delta,\delta]\subset\delta$ but $\D^1\h\nsubseteqq\delta$, so that $[\delta,Z]$ is not included in $\delta$. Therefore there exists $X\in\delta$ such that $\ad(X)Z=aZ$ $(\mbox{mod }{\delta})$, for some real number $a\neq 0$. We can choose $a=1$ so that
$$
\forall k\geq 1 \qquad \ad^k(X)Z=Z \ \  (\mbox{mod }{\delta})
$$
which proves that $\h$ is not nilpotent, a contradiction.
\item Let now $\g$ and $\Delta$ as in the statement. If $\Delta$ is an ideal of $\g$ then the result comes from Proposition \ref{Omegaclosed}. If not $\D^1\g$ is not included in $\Delta$ and there exists $Z\in\D^1\g$ such that $\g=\R Z\oplus \Delta$. Moreover the rank condition implies that $D(\Delta)\nsubseteqq \Delta$.
\item Let $\eta$ be the largest nilpotent ideal of $\g$ (see \cite{Bourbaki1}, $\S\ 4$, $n^o\ 4$). It is known that for any derivation $d$ of the solvable algebra $\g$, the inclusion $d(\g)\subset \eta$ holds, so that $\eta$ is a characteristic ideal of $\g$ (see \cite{Bourbaki1}, $\S\ 5$, Proposition 6 and its corollary). Moreover $\D^1\g$ being a nilpotent ideal of $\g$ is included in $\eta$, hence $Z\in \eta$ and $\g=\eta+\Delta$. Consequently
$$
\dim (\eta\ \cap\Delta)=\dim(\eta)+\dim(\Delta)-\dim(\g)=\dim(\eta)-1,
$$
which proves that $\eta\ \cap\Delta$ is a codimension $1$ subalgebra of $\eta$. Since $\eta$ is nilpotent, and according to Item 1, $\eta\ \cap\Delta$ is an ideal of $\eta$. Actually $\eta\ \cap\Delta$ turns out to be an ideal of $\g$. Indeed any $Y\in \g$ writes $Y=X+aZ$ where $X\in \Delta$ and $a\in\R$. Since $Z\in \eta$ and $\eta\ \cap\Delta$ is an ideal of $\eta$ it is clear that $\ad(Z)$ sends $\eta\ \cap\Delta$ into $\eta\ \cap\Delta$. The same is true for $\ad(X)$ because $X$ belongs to the algebra $\Delta$ and $\eta$ is an ideal of $\g$. Finally $\eta\ \cap\Delta$ is invariant for all inner derivations, hence an ideal of $\g$.
\item Let us show that $D^{-1}\Delta$ is a subalgebra. Firstly we have $D(D^{-1}\Delta)\subset D(\g)\subset \eta$ and $D(D^{-1}\Delta)\subset \Delta$ by definition, so that $D(D^{-1}\Delta)\subset \eta\ \cap\Delta$. Then for all $X_1$, $X_2$ in $D^{-1}\Delta$ we get $DX_i\in \eta\ \cap\Delta$, hence $[DX_i,X_j]\in \eta\ \cap\Delta$ for $i,j=1,2$. Consequently
$$
D[X_1,X_2]=[DX_1,X_2]+[X_1,DX_2] \in \eta\ \cap\Delta
$$
and $[X_1,X_2]\in D^{-1}\Delta$, which finishes the proof.
\item In this item we show that $F(g)\in\eta$ for all $g\in G$. According to Proposition \ref{A2} of the appendix
$$
\forall Y\in \g, \quad F(\exp tY)=\sum_{k\geq 1}(-1)^{k-1}\frac{t^k}{k!}\ad^{(k-1)}(Y)DY\ \in \eta
$$
because $DY$ belongs to the ideal $\eta$. Secondly
$$
F(g\exp tY)=F(\exp tY)+e^{-t\ad(Y)}F(g)
$$
belongs to $\eta$ as soon as $F(\exp tY)$ and $F(g)$ both belong to $\eta$. Since $g=\exp(Y_1)\dots\exp(Y_k)$ for some $Y_1,\dots,Y_k\in \g$ the result is obtained by induction.
\item To finish the proof of the theorem let $g$, $g'$ in $\Z$. On the one hand $F(g)$ belongs to the ideal $\eta\ \cap\Delta$. On the other one $g'=\exp(X_1)\dots\exp(X_k)$ for some $X_1,\dots,X_k\in \g$, hence
$$
\Ad(g')=e^{\ad(X_1)}\dots e^{\ad(X_k)}
$$
and $\Ad(g')F(g)$ belongs to $\eta\ \cap\Delta$. The singular locus $\Z$ is a subgroup of $G$ according to Lemma \ref{lemmeAd}.
\end{enumerate}

\hfill $\Box$

\vskip 0.2cm

\noindent{\bf Remarks}.
\begin{enumerate}
\item  The beginning of the proof shows that when $\g$ is nilpotent then $\Delta$ is an ideal as soon as it is a subalgebra. In the solvable case it is easy to exhibit examples where $\Delta$ is a subalgebra but not an ideal (see Example \ref{ZAff} on $Aff_+(2)$).
\item The same is true on semi-simple Lie groups since a codimension one subalgebra cannot be an ideal. But apart from the fact that codimension one subalgebras are not very common in semi-simple Lie algebras, that condition is not here sufficient for $\Z$ to be a group (see Example \ref{ZSL2} on $SL2$).
\end{enumerate} 

\vskip 0.2cm

There is another case where we can assert that $\Z$ is a subgroup, it is when it is equal to the set of singularities of $\X$.

\newtheorem{ZequalChi}[Omegaclosed]{Proposition}
\begin{ZequalChi} \label{ZequalChi} 
If $\Z$ is equal to the set $\Z_{\X}$ of singularities of $\X$ then it is a closed Lie subgroup of $G$ whose Lie algebra is equal to $\ker(D)$.
\end{ZequalChi}

\demo
The set $\Z_{\X}$ is obviously a closed subgroup of $G$, since
$$
\Z_{\X}=\bigcap_{t\in\R}\{g\in G;\ \varphi_t(g)=g\},
$$
where $(\varphi_t)_{t\in\R}$ stands for the flow of $\X$. Moreover
$$
\begin{array}{l}
\forall s\in\R \ \ \exp(sY)\in \Z_{\X}\\
\Longleftrightarrow \forall s\in\R, \  \forall t\in\R,\  \exp(sY)=\varphi_t(\exp(sY))=\exp(e^{tD}sY)\\
\Longleftrightarrow \forall t\in\R,\ Y=e^{tD}Y\\
\Longleftrightarrow DY=0,
\end{array} 
$$
which shows that the Lie algebra of $\Z_{\X}$ is $\ker D$.

\hfill $\Box$

\vskip 0.2cm

This case is illustrated in Examples \ref{ZAff} and \ref{ZHeisenberg}.


\subsection{Necessary conditions for $\Z$ to be a subgroup}
When the one form $D^*\omega$ does not vanish, the singular locus is locally, around the identity, a submanifold, whose tangent space at $e$ is $T_e\Z=\ker(D^*\omega)=D^{-1}\Delta$. One could expect  $\Z$ to be a subgroup if $D^{-1}\Delta$ is a subalgebra, but this statement is wrong, even locally, and even if $D^{-1}\Delta$ is an ideal of $\g$ (see the examples of Section \ref{SLE}, in particular Example 3.5.3.1).

In order to state a necessary condition for $\Z$ to be a subgroup, and a locally sufficient one, we will use the following stronger algebraic condition:

Let $\h$ be a subalgebra of $\g$. It will be said to satisfy the \textbf{condition $(HZ)$} if
$$
\begin{array}{l}
\forall\ Y\in \h, \ \forall\ m\geq 0,\ \forall\  Z_1,\dots,Z_m\in \h,
\qquad\langle\omega, \ad(Z_1)\dots\ad(Z_m)DY\rangle=0.
\end{array}
$$

\newtheorem{Analyse}[Omegaclosed]{Proposition} 
\begin{Analyse} \label{Analyse}
If $\Z$ is a group, then it is a Lie subgroup of $G$ whose Lie algebra $\z$ satisfies Condition $(HZ)$.
\end{Analyse}

\demo The set $\Z$ being closed is a Lie subgroup as soon as it is a group.

For all $Z_1,\dots,Z_m\in\z$ and for all $t_1,\dots,t_m\in \R$ the point $\exp(t_1Z_1)\dots\exp(t_mZ_m)$ belongs to $\Z$. Since $\Ad(\exp(t_1Z_1)\dots\exp(t_mZ_m))=e^{t_1\ad(Z_1)}\dots e^{t_m\ad(Z_m)}$, and according to Lemma \ref{lemmeAd}, we get for all $g\in\Z$
$$
\langle\omega,e^{t_1\ad(Z_1)}\dots e^{t_m\ad(Z_m)}F(g)\rangle=0.
$$
Derivating this equality at $t_1=0,\dots,t_m=0$, and thanks to the linearity of $\ad(Z_i)$, we obtain
$$
\langle\omega,ad(Z_1)\dots \ad(Z_m)F(g)\rangle=0.
$$
To finish let $g=\exp(tY)$, where $Y\in\z$. According to Formula (\ref{TeF}) (Proposition \ref{A1}, Appendix), we have
$$
\dto F(\exp tY)=DY,
$$
which gives the desired equality:
$$
\langle\omega,ad(Z_1)\dots \ad(Z_m)DY\rangle=0.
$$

\hfill $\Box$

\vskip 0.3cm

\noindent{\bf Remarks}
\begin{enumerate}
\item Proposition \ref{Analyse} is a local result: actually it holds as soon as $\Z$ is locally, in a neighbourhood of the identity, equal to some Lie subgroup of $G$. 
\item The successive derivatives of $F(\exp(tY))$ do not provide more information since some of the $Z_i$'s can be chosen equal to $Y$.
\item Condition (HZ) is trivially satisfied only if $\Delta$ is an ideal.
\end{enumerate}

In order to investigate the converse to Proposition \ref{Analyse} we will consider a subalgebra $\h$ of $\g$, assumed to satisfy Condition (HZ), and we will associate to $\h$ the following subspace of $\Delta$:
$$
\delta_h=\{X\in\Delta;\  \forall\ m\geq 0,\ \forall\  Z_1,\dots,Z_m\in \h
\quad\ \ad(Z_1)\dots\ad(Z_m)X\in \Delta\}.
$$
Notice that $\delta_h$ depends on $\h$ and, thanks to Condition (HZ), that $D\h\subset \delta_h$. It is also clear that $\delta_h$ is $\ad(Z)$-invariant for all $Z\in \h$.

\newtheorem{Synthese}[decompositionDez]{Lemma} 
\begin{Synthese} \label{Synthese}
Let $\h$ be a subalgebra of $\g$ that satisfies Condition (HZ), and let $H$ be the connected subgroup generated by $\h$. Then
\begin{enumerate}
\item The Lie subgroup $H$ is included in $\Z$.
\item Let $g\in\Z$. If $F(g)\in \delta_h$, then $gH$ is included in $\Z$.
\end{enumerate}
\end{Synthese}

\demo
\begin{enumerate}
\item Let $Y\in\h$ and $k\geq 1$. Condition (HZ) with $m=k-1$ and $Z_i=Y$ for $i=1,\dots,k-1$ is:
$$
\ad^{(k-1)}(Y)DY\in \delta_h\subset \Delta
$$
Consequently for all $t\in \R$
$$
F(\exp tY)=\sum_{k\geq 1}(-1)^{k-1}\frac{t^k}{k!}\ad^{(k-1)}(Y)DY\in \delta_h\subset \Delta
$$
and $\exp(tY)$ belongs to $\Z$ for all $t\in\R$.
\item Now let $g\in \Z$ such that $F(g)\in \delta_h$. For any $Y\in \h$ we have
$$
\begin{array}{ll}
\langle\omega, F(g\exp tY)\rangle & = \langle\omega, F(\exp tY)+e^{-t\ad(Y)}F(g)\rangle\\
& =\langle\omega, F(\exp tY)\rangle+\sum_{k\geq 0}\frac{(-t)^k}{k!}\langle\omega,\ad^{(k)}(Y)F(g)\rangle\\
& =0
\end{array}
$$
because on the one hand the first term vanishes according to the first item, and on the other one $\ad^{(k)}(Y)F(g)$ belongs to $\Delta$ since $F(g)$ belongs to $\delta_h$. Notice moreover that $F(g\exp tY)$ belongs to $\delta_h$, since $\delta_h$ is $\ad(Z)$-invariant for all $Z\in \h$. 
\item
Let $Z_1,\dots,Z_m$ in $\h$ and $t_1,\dots,t_m$ in $\R$. According to Item 2 we obtain by induction that $\exp(t_kZ_k)\dots\exp(t_1Z_1)$ belongs to $\Z$ for $k=1,\dots,m$, but moreover that\\ $F(\exp(t_kZ_k)\dots\exp(t_1Z_1))$ belongs to $\delta_h$, allowing the induction.

This shows that $H$ is included in $\Z$.
\item To finish let $g\in\Z$ such that $F(g)\in \delta_h$. We get $gH\subset \Z$ by the same reasoning as in Item 3.
\end{enumerate}

\hfill $\Box$

\newtheorem{localsubgroup}{Corollary}
\begin{localsubgroup}\label{localsubgroup}
If $D^*\omega$ is not zero and if $\z=\ker(D^*\omega)$ is a subalgebra of $\g$ that satisfies Condition (HZ), then the connected subgroup $\mbox{Gr}(\z)$ generated by $\z$ is included in $\Z$.

Moreover there exists a neighbourhood $V$ of the identity such that:
$$
\Z\cap V=\mbox{Gr}(\z)\cap V.
$$
\end{localsubgroup}

\demo

The first part is immediate from Lemma \ref{Synthese}. For the second one the condition $D^*\omega\neq 0$ implies that $\Z$ is locally, in a neighbourhood $V$ of $e$, a $(n-1)$-dimensional submanifold (recall that $\Z=\{\psi=0\}$ and  $T_e\psi=D^*\omega$). This submanifold containing the $(n-1)$-dimensional Lie group 
$\mbox{Gr}(\z)$, is equal to that last within $V$.

\hfill $\Box$

The examples 3.5.3.1 on the group Heisenberg (see Section \ref{ZHeisenberg}) show that the assumptions of Corollary \ref{localsubgroup} are not sufficient to know what happens far from the identity: they imply neither that $\psi$ is regular in $\Z$ nor that one has $F(g)\in\delta_{\z}$ for all the points $g\in\Z$, even in the connected component of $\Z$.


\subsection{Examples of singular loci}\label{SLE}
\subsubsection{Abelian groups}
Consider the ARS on $\R^n$ defined by $(A,b_{1},...,b_{n-1} )$ where $A$ is a real matrix of
order $n$ and  $b_{1},...,b_{n-1}$ are
$n-1$ linearly independent constant vector fields . It is straightforward to show that the singular locus is a codimension one subspace of $\R^n$.


\subsubsection{The 2D affine group}\label{ZAff}
 Let $G$ be the connected component of $e$ in the $2$-dimensional affine group:
 $$
 G=Aff_+(2)=\left\{\begin{pmatrix}x&y\\0&1\end{pmatrix};\ \ (x,y)\in\R_+^*\times \R\right\}.
 $$
 Its Lie algebra $\g=\mathfrak{aff}(2)$ is solvable, generated by
 $$
 X=\begin{pmatrix}1&0\\0&0\end{pmatrix} \quad \mbox{ and }\quad Y=\begin{pmatrix}0&1\\0&0\end{pmatrix}
 $$
 with $[X,Y]=XY-YX=Y$, and identified with the set of left-invariant vector fields which is therefore generated by
 $$
 gX=\begin{pmatrix}x&0\\0&0\end{pmatrix} \quad \mbox{ and }\quad gY=\begin{pmatrix}0&x\\0&0\end{pmatrix}\quad \mbox{ where} \quad g=\begin{pmatrix}x&y\\0&1\end{pmatrix}.
 $$

 All the derivations are inner and the one whose matrix in the basis $(X,Y)$ is
$\begin{pmatrix}
0 & 0 \\ 
a & b%
\end{pmatrix}$ is equal to $D=-\ad(aY-bX)$. The associated linear vector field $\X$ is given at the point $g$ by 
\begin{equation*}
\X_{g}=%
\begin{pmatrix}
0 & a(x-1)+by \\ 
0 & 0%
\end{pmatrix}%
.
\end{equation*}%
Let us consider the simple ARS defined by $\X$ and $Y_1=\alpha X+\beta Y=\begin{pmatrix} \alpha x & \beta x \\ 0 & 0\end{pmatrix}$. In order that $\X$ and $Y_1$ be independant on some subset the constant $\alpha$ must not vanish.

Then $DY_1=a\alpha Y+b\beta Y$ and the rank condition
is equivalent to $\alpha a+\beta b\neq 0$.

As $x>0$ the first component of $Y_1$, that is $\alpha x$, never vanishes, so that $Y_1$ and $\X$ are colinear if and only if $\X$ vanishes. Finally we get
$$
\Z=\left\{g\in G;\ \X_{g}=0\right\}=\{(x,y)\in\R_+^*\times \R;\ a(x-1)+by=0\} .
$$
According to Proposition \ref{ZequalChi} it is a closed Lie subgroup of G. Another proof comes from the fact that $\Delta$ is always a subalgebra since it is one dimensional. As the group $Aff_+(2)$ is solvable we know by Theorem \ref{Zgroupe} that $\Z$ should be a subgroup.


\subsubsection{The Heisenberg group}\label{ZHeisenberg}

Let $G$ be the simply connected Heisenberg Lie group of
dimension three 
\begin{equation*}
G=\left\{ 
\begin{pmatrix}
1 & x & z \\ 
0 & 1 & y \\ 
0 & 0 & 1%
\end{pmatrix}%
;\ x,y,z\in 
\R
\right\} \text{.}
\end{equation*}

Its Lie Algebra $\g$ is generated by 
\begin{equation*}
X=
\begin{pmatrix}
0 & 1 & 0 \\ 
0 & 0 & 0 \\ 
0 & 0 & 0%
\end{pmatrix}%
,\;Y=
\begin{pmatrix}
0 & 0 & 0 \\ 
0 & 0 & 1 \\ 
0 & 0 & 0%
\end{pmatrix}%
,\text{ }Z=
\begin{pmatrix}
0 & 0 & 1 \\ 
0 & 0 & 0 \\ 
0 & 0 & 0%
\end{pmatrix}%
\end{equation*}
where $[X,Y]=XY-YX=Z$ and the other brackets vanish.
As left-invariant vector fields they write in natural coordinates:
\begin{equation*}
X=\frac{\partial }{\partial x},\qquad Y=\frac{\partial }{\partial y}+x\frac{%
\partial }{\partial z},\qquad Z=\frac{\partial }{\partial z}.
\end{equation*}%
The derivations of $\g$ are the endomorphisms $D$ whose matrix in the basis $(X,Y,Z)$ has the form:
\begin{equation*}
D=
\begin{pmatrix}
a & b & 0 \\ 
c & d & 0 \\ 
e & f & a+d%
\end{pmatrix}%
\end{equation*}%
and the associated linear vector
field is:
\begin{equation*}
\X(g)=(ax+by)\frac{\partial }{\partial x}+(cx+dy)\frac{\partial }{\partial y}%
+(ex+fy+(a+d)z+\frac{1}{2}cx^{2}+\frac{1}{2}by^{2})\frac{\partial }{%
\partial z}.
\end{equation*}%
More details can be found in \cite{Jou09}.

Let us now consider a simple ARS on the group Heisenberg, defined by a $2$-dimensional subspace $\Delta$ of $\g$ and a derivation $D$. We get two very different behaviour depending on whether $\Delta$ is a subalgebra or not.

One can also notice that the number of connected components of $G\setminus \Z$ ranges from 1 to 4.
\begin{enumerate}
\item \textbf{$\Delta$ is a subalgebra}. According to Theorem \ref{Zgroupe} the singular locus $\Z$ is a subgroup of $G$. Up to an automorphism of $\g$ we can assume that $\Delta=\Span\{X,Z\}$. Then the rank condition is satisfied if and only if $c\neq 0$ (in the matrix of $D$). A particular case is when $\Z$ is exactly the set of singularities of $\X$. One example is obtained by
$$
D=
\begin{pmatrix}
0 & 0 & 0 \\ 
1 & 0 & 0 \\ 
e & 0 & 0%
\end{pmatrix}
$$
but is far from being the only one. Here $\X=x\frac{\partial }{\partial y}
+(ex+\frac{1}{2}x^{2})\frac{\partial }{\partial z}$ and it is easily seen that $\Z=\{\X=0\}=\{x=0\}$.
\item \textbf{$\Delta$ is not a subalgebra}. We can assume without lost of generality that $\Delta=\Span\{X,Y\}$. It is easy to see that the ARS is well defined for any derivation different from zero. Indeed the rank condition is satisfied and
$$
\Z=\{ex+fy+(a+d)z-\frac{1}{2}cx^{2}+\frac{1}{2}by^{2}-dxy=0\}
$$
is equal to $\R^3$ if and only if all the coefficients of $D$ vanish.

It is also clear that the loci defined by these quadratic forms need not be subgroups, not even submanifolds. Let us exhibit some particular examples.
\vskip 0.2cm
\noindent \textbf{Example 3.5.3.1}
Let $D$ have the following form:
$$
D=\begin{pmatrix}
a&b&0\\c&-a&0\\0&1&0
\end{pmatrix}.
$$

In the $\{X,Y,Z\}$ basis we have $\omega=(0,0,1)$ and $D^*\omega=(0,1,0)$.

The one form $D^*\omega$ does not vanish, and its kernel, which is also the tangent space to $\Z$ at the origin, is the subalgebra $\z=D^{-1}\Delta=\Span\{X,Z\}$. However this fact is not sufficient for $\Z$ to be a subgroup, even locally.

If $c\neq 0$ the condition $(HZ)$ is not satisfied, and $\Z=\{y-\frac{1}{2}cx^{2}+\frac{1}{2}by^{2}+axy=0\}$ is not a local subgroup around the origin: it is straightforward to check that $(x,y,z)^{-1}=(-x,-y,-z+xy)$ does not belong to $\Z$ in general when $(x,y,z)\in\Z$.

 If $c=0$ then $\z$ satisfies Condition $(HZ)$ so that we know by Corollary \ref{localsubgroup} that $\Z$ contains the group generated by $\z$, that is $\mbox{Gr}(\z)=\{y=0\}$. However
$$
\Z=\{y(1+\frac{1}{2}by+ax)=0\}
$$
and $\Z$ reduces to $\{y=0\}$ if and only if $a=b=0$. Otherwise $\Z$ is the union of the plane $\{y=0\}$ with another plane which may intersect it (if $a\neq 0$) or be parallel to it (if $a=0$).

Let us consider the case $a=1$.

At the points $g=(-1,0,z)$, that belong to $\mbox{Gr}(\z)$, the function $\psi$ is singular, that is $T_g\psi$ vanishes.

At the points $g=(-1-\frac{1}{2}by,y,z)$ with $y\neq 0$, it is clear that $g\mbox{Gr}(\z)$, the translation by $g$ of the connected group generated by $\z$, is not included in $\Z$.

\vskip 0.2cm
\noindent \textbf{Example 3.5.3.2. Degenerated case}

Consider
$$
D=\begin{pmatrix}
0&b&0\\c&0&0\\0&0&0
\end{pmatrix} \quad\mbox{ with } b>0 \mbox{ and } c<0
$$
Here $D^{-1}\Delta=\g$ and the condition $(HZ)$ is not satisfied.

If $b>0$ and $c<0$, then $\Z$ is the codimension $2$ subgroup $\{x=y=0\}$ of $\g$.

If $b$ and $c$ are both positive (or negative), then $\Z=\{bx=\pm by\}$ is the union of two secant planes.

If $b\neq 0$ and $c=0$, then $\Z$ is the codimension $1$ subgroup $\{y=0\}$ of $\g$.
\vskip 0.2cm
\noindent \textbf{Example 3.5.3.3. Tangential case}

Let
$$
D=\begin{pmatrix}
a&b&0\\c&d&0\\0&0&a+d
\end{pmatrix} \quad\mbox{ with } a+d\neq 0.
$$
The one forms $\omega=\begin{pmatrix}0 & 0 & 1\end{pmatrix}$ and $D^*\omega=\begin{pmatrix}0 & 0 & a+d\end{pmatrix}$ define the same subspace $\Delta=D^{-1}\Delta$ of $\g$. As the differential at the origin of the defining function of $\Z$, that is $T_e\Psi=D^*\omega$ does not vanish, the singular locus is a submanifold in a neighbourhood of $e$, and $e$ is a tangential point, that is $\Delta=T_e\Z$.

Unlike the generic case (see \cite{BCGM14}) the tangential points need not be isolated. Consider
$$
D=\begin{pmatrix}
0&0&0\\2&1&0\\0&0&1
\end{pmatrix}.
$$
The singular locus is $\Z=\{z=x^2+xy\}$, and $\Delta$ is tangent to $\Z$ along the parabola $\{z=-x^2\}$ contained in the plane $\{y=-2x\}$. 
\end{enumerate}


\subsubsection{The special linear group $SL(2;\R)$}\label{ZSL2}
\noindent \textbf{Example 3.5.4.1}. 
Let $G=SL(2;\mathbb{R})$ be the order $2$ special linear group. Its Lie algebra $%
\mathfrak{g}=\mathfrak{sl}(2;\mathbb{R})$ is the set of matrices with trace zero. The usual basis:
\begin{equation*}
H=\begin{pmatrix} 1 & 0 \\  0 & -1 \end{pmatrix},
\quad
X=\begin{pmatrix} 0 & 1 \\  0 & 0 \end{pmatrix},
\quad
Y=\begin{pmatrix} 0 & 0 \\ 1 & 0 \end{pmatrix}
\end{equation*}
of $\mathfrak{sl}(2;\mathbb{R})$ verifies
\begin{equation*}
\left[ H,X\right] =2X,\quad \left[ H,Y\right] =-2Y,\quad
\left[ X,Y\right] =H.
\end{equation*}
Consider the simple ARS defined by $\Delta=\Span\{H,X\}$ and the derivation $D=-\ad(Y)$ ($\mathfrak{sl}(2;\mathbb{R})$ being semi-simple, all the derivations are inner). We get $D^{-1}\Delta=\Span\{X,Y\}$ which is not a subalgebra. Consequently the singular locus cannot be a codimension one subgroup of $SL(2;\R)$ despite the fact that $\Delta$ is a subalgebra of $\mathfrak{sl}(2;\mathbb{R})$. This shows that Theorem \ref{Zgroupe} does not apply to general groups.

Let us compute the singular locus.
At the point $g=\begin{pmatrix} a & b \\  c & d \end{pmatrix}$ with $ad-bc=1$, the linear vector field associated to $D=-\ad(Y)$ is $\X_{g}=gY-Yg=\begin{pmatrix} b & 0 \\ d-a & -b \end{pmatrix}$. A straightforward computation shows that the vectors $\X_g,H_g,X_g$ are linearly dependent if and only if $a=\pm 1$, in other words that
\begin{equation*}
\Z=\left\{ a=\pm 1\right\}.
\end{equation*}%
Thus the singular locus $\Z$ is a submanifold of $SL(2;\R)$ but nor a subgroup neither a connected set.

Notice that the ARS is well defined: indeed $\Z\neq G$ and the rank condition is satisfied.

\vskip 0.2cm

\noindent \textbf{Example 3.5.4.2}. 
Another example can be obtained by switching the roles of $\Delta$ and $\D^{-1}\Delta$. More accurately let $\Delta=\Span\{X,Y\}$, it is not a subalgebra, and $D=-\ad(X)$. We get easily $D^{-1}\Delta=\Span\{H,X\}$. It is a subalgebra of $\mathfrak{sl}(2;\mathbb{R})$ that satisfies Condition $(HZ)$, hence the singular locus contains the connected subgroup generated by $D^{-1}\Delta$. However this locus turns out to be $\Z=\{cd=0\}$, which has three connected components:
$$
C_1=\left\{\begin{pmatrix}
a & b \\ 0 & \frac{1}{a}
\end{pmatrix}
;\ a>0\right\},
\qquad
C_2=\left\{\begin{pmatrix}
a & b \\ 0 & \frac{1}{a}
\end{pmatrix}
;\ a<0\right\},
\quad
$$
$$
C_3=\left\{\begin{pmatrix}
a & b \\ c & 0
\end{pmatrix}
;\ bc=-1\right\}.
$$
$C_1$ is a connected subgroup, $C_1\cup C_2$ is a non connected subgroup, but $\Z=C_1\cup C_2\cup C_3$ is clearly not a subgroup of $SL(2;\R)$.


\section{Hamiltonian equations}\label{HE}

The aim of this section is to state the Hamiltonian equations of the PMP applied to ARS's on Lie groups. We follow the same lines as in the invariant case (see for instance \cite{AS04}) but here the Hamiltonian depends on the point $g\in G$  which entails some more complications (see also \cite{Sachkov09}).


\subsection{The canonical symplectic structure of $T^*G$}
All the material of this subsection is standard and can be found in the previous references.

As usual the cotangent bundle $T^*G$ is identified to $\g^*\times G$ by
$$
\begin{array}{lll}
\Phi\ : & \g^*\times G & \longrightarrow T^*G\\
        & (\lambda,g) & \longmapsto \bar{\lambda}_g=\lambda\circ TL_{g^{-1}} \in T_g^*G,
\end{array}
$$
so that the equality $\ds \left\langle \bar{\lambda}_g,Y_g \right\rangle=\left\langle\lambda,Y \right\rangle$ holds for all $Y \in \g$ (recall that $Y_g$ stands for $TL_g.Y$).

The projection from $T^*G$ to $G$ is denoted by $\Pi$ in what follows.

The tautological one-form $s$ on $T^*G$ being defined by $\ds \left\langle s_{\bar{\lambda}},\zeta\right\rangle=\left\langle \bar{\lambda}, \Pi_*\zeta\right\rangle$, its pullback by $\Phi$ is the one-form $\Phi^*s$ defined on $\g^*\times G$ by (here $X\in\g$):
$$
\left\langle \Phi^*s_{(\lambda,g)},(\xi,X_g)\right\rangle=\left\langle \lambda,X\right\rangle.
$$
In this setting, the symplectic form $d\Phi^*s=\Phi^*ds$ is characterized by
$$
\Phi^*ds_{(\lambda,g)}\left((\xi,X_g),(\eta,Y_g)\right)=\left\langle \xi,Y\right\rangle-\left\langle \eta,X\right\rangle-\left\langle \lambda,[X,Y]\right\rangle.
$$

To finish this review let $h\in \cci(T^*G)$ and let $\vec{h}$ be the Hamiltonian vector field associated to $h$ by $\ds dh=-i_{\vec{h}}ds$. If $H$ is equal to $h\circ \Phi$ and $\vec{H}$ is defined by $\vec{h}=\Phi_*\vec{H}$, it is easily verified that $\vec{H}$ is the Hamiltonian vector field associated to $H$, i.e.
$$
\ds dH=-i_{\vec{H}}\Phi^*ds.
$$


\subsection{Computation of the Hamiltonian vector fields}

In what follows we set $\ds \sigma=\Phi^*ds$.

Let $H$ be an Hamiltonian on $\g^*\times G$ and $\vec{H}$ the associated Hamiltonian vector field. We can identify the tangent space to $\g^*\times G$ at the point $(\lambda,g)$ with $\g^*\times T_gG$, and write $\vec{H}=(\xi,X_g)$ at this point. Then for all $(\eta,Y_g)\in \g^*\times T_gG$ we get on the one hand
$$
\ds dH(\eta,Y_g)=\left\langle \frac{\partial H}{\partial \lambda}, \eta\right\rangle+\left\langle \frac{\partial H}{\partial g},Y_g\right\rangle,
$$
and on the other hand
$$
\begin{array}{ll}
\ds dH(\eta,Y_g) & =-\sigma_{(\lambda,g)}\left((\xi,X_g)(\eta,Y_g)\right)\\
\ds              & =-\left\langle \xi,Y\right\rangle+\left\langle \eta,X\right\rangle+\left\langle \lambda,[X,Y]\right\rangle.
\end{array}
$$
\begin{enumerate}
	\item Setting $Y=0$ we get $\ds X=\frac{\partial H}{\partial \lambda}(\lambda,g)\in (\g^*)^*\equiv \g$.
	\item Setting $\eta=0$ we get
	$$
	\left\langle \frac{\partial H}{\partial g},Y_g\right\rangle=-\left\langle \xi,Y\right\rangle+\left\langle\lambda,[X,Y]\right\rangle
	$$
	hence
	$$
	\left\langle \xi,Y\right\rangle=\left\langle(\ad(X))^*\lambda,Y\right\rangle-\left\langle (TL_g)^*\frac{\partial H}{\partial g},Y\right\rangle
	$$
	so that
	$$
	\xi=(\ad(X))^*\lambda-(TL_g)^*\frac{\partial H}{\partial g}.
	$$
\end{enumerate}
Summarizing we obtain
$$
\ds \vec{H}=\left(\frac{\partial H}{\partial \lambda}\right)_g\frac{\partial }{\partial g}
+\left(\left(\ad (\frac{\partial H}{\partial \lambda})\right)^*\lambda-(TL_g)^*\frac{\partial H}{\partial g}\right)\frac{\partial}{\partial \lambda},
$$
where $\ds \left(\frac{\partial H}{\partial \lambda}\right)_g$ should be understood as $\ds TL_g\frac{\partial H}{\partial \lambda}$.


\subsection{Left-invariant and linear vector fields} \label{lilvf}

The Hamiltonian associated to a left-invariant vector field $Y$, that is $H=\left\langle \lambda,Y\right\rangle$, does not depend on $g$, so that the corresponding Hamiltonian equations turn out to be
$$
\left\{
\begin{array}{ll}
\ds \dot{g} & =Y_g\\
\ds \dot{\lambda} & =(\ad(Y))^*\lambda
\end{array}
\right.
$$
Let us now consider a linear vector field $\X$, whose associated derivation is $D=-\ad(\X)$, and let us define the Hamiltonian
$$
H(\lambda,g)=\left\langle \lambda_g,\X_g\right\rangle=\left\langle \lambda,TL_{g^{-1}}\X_g\right\rangle=
\left\langle \lambda,F_g\right\rangle
$$
According to Corollary \ref{C1} (in the Appendix) we have for any $Y\in\g$
$$
\ds \frac{\partial H}{\partial g}.Y_g
=\left\langle \lambda, DY+\ad(F_g)Y\right\rangle
 =\left\langle \left(D+\ad(F_g)\right)^*\lambda,Y\right\rangle.
$$

Finally the Hamiltonian equations for the Hamiltonian $H=\left\langle \lambda,F(g)\right\rangle$ associated to the linear vector field $\X$ are:

$$
\left\{
\begin{array}{ll}
\ds \dot{g} & =\X_g\\
\ds \dot{\lambda} & =\left(D+\ad(F_g)\right)^*\lambda
\end{array}
\right.
$$

\vskip 0.2cm

\noindent \textbf{Remark: the inner case}

If the derivation $D$ is inner, that is if $D=-\ad X$ for some $X\in \g$, then $\X$ has the form $\X_g=TL_g.X-TR_g.X$ and $F(g)=X-\Ad(g^{-1})X$. Consequently
$$
D+\ad (F_g)=D+\ad X-\ad(\Ad(g^{-1})X)=-\ad(\Ad(g^{-1})X),
$$
and the second equation reduces to
$$
\ds \dot{\lambda}=\left(-\ad(\Ad(g^{-1})X)\right)^*\lambda.
$$


\subsection{Hamiltonian equations of a simple ARS}

Consider an ARS defined as previously by $(\X,Y_1,\dots,Y_{n-1})$, and consider the Hamiltonian
$$
\ds \mathcal{H}_{\nu}(\lambda,g,v,u_1,\dots,u_{n-1})=\left\langle \lambda,v\X+\sum_{1}^{n-1}u_jY_j)\right\rangle
\ds -\frac{1}{2}\nu\left(v^2+\sum_{1}^{n-1}u_j^2\right).
$$
The associated equations are
$$
\left\{
\begin{array}{ll}
\ds \dot{g} & =v\X+\sum_{1}^{n-1}u_jY_j\\
\ds \dot{\lambda} & =\left(vD+\ad\left(vF(g)+\sum_{1}^{n-1}u_jY_j\right)\right)^*\lambda
\end{array}
\right.
$$
In particular the equations of the normal extremals are obtained by application of the Pontryagyn Maximum Principle (see for instance \cite{PBGM62}, \cite{ABB14}, \cite{AS04}, \cite{Jurdjevic97}) with $\nu=1$. As usual the maximization of $\mathcal{H}_1$ w.r.t. $v,u_1,\dots,u_{n-1}$ gives $v=\langle\lambda,\X \rangle$ and $u_j=\langle\lambda,Y_j \rangle$ for $j=1,\dots,n-1$. The maximized Hamiltonian is:
$$
\ds H_1(\lambda,g)=\frac{1}{2}\left\langle \lambda,\X\right\rangle^2+\frac{1}{2}\sum_{1}^{n-1}\left\langle \lambda,  Y_j\right\rangle^2.
$$

The case of the abnormal extremals is treated in Section \ref{Abnormal}.


\subsection{Semi-simple Lie groups}

In the semi-simple case all the derivations are inner and the remark of the end of Section \ref{lilvf} applies. Moreover there exist on $\g$ an invariant scalar product. In order to avoid confusion it will be denoted by $\ds \left\langle .\ ,. \right\rangle_s$ while the duality bracket will be denoted by $\left\langle .\ ,. \right\rangle_d$, and the invariance of the scalar product means
$$
\forall X,Y,Z\in \g \qquad \left\langle \ad(X)Y,Z\right\rangle_s=-\left\langle Y,\ad(X)Z \right\rangle_s.
$$
Thanks to this scalar product, we can identify $\g$ and $\g^*$ by $X\in\g\longmapsto \lambda_X=\left\langle X,.\right\rangle_s\in\g^*$, so that
$$
\begin{array}{ll}
\left\langle \ad(X)^*\lambda_Y,Z\right\rangle_d & = \left\langle \lambda_Y,\ad(X)Z\right\rangle_d\\
 & = \left\langle Y,\ad(X)Z\right\rangle_s\\
 & = -\left\langle \ad(X)Y,Z\right\rangle_s
\end{array}
$$
and $\ad(X)^*$ can be identified with $-\ad(X)$. The previous equations become:
$$
\left\{
\begin{array}{ll}
\ds \dot{g} & =v\X+\sum_{1}^{n-1}u_jY_j\\
\ds \dot{Z} & =[v\Ad(g^{-1})X-\sum_{1}^{n-1}u_jY_j,Z]
\end{array}
\right.
$$
where $Z\in \g$.


\section{Abnormal extremals}\label{Abnormal}

Let $(\lambda(t),g(t))$, $t\in[0,T]$, be an abnormal extremal. We know that
\begin{enumerate}
\item[(i)] $\lambda(t)$ does not vanish;
\item[(ii)] $\langle\lambda(t),Y_i\rangle\equiv 0$ for $i=1,\dots,n-1$;
\item[(iii)] $\langle\lambda(t)\circ TL_{g(t)^{-1}},\X_{g(t)}\rangle=\langle\lambda(t),F_{g(t)}\rangle\equiv 0$.
\end{enumerate}
These conditions imply that $g(t)$ belongs to the singular set $\Z$ for all $t$. Consequently the linear vector field $\X$ is a linear combination of $Y_1,\dots,Y_{n-1}$ along the curve $g(t)$ and we can assume the control $v$ to vanish. This may modify the optimality of the associated control, but not the geometry of the curve, which is what we are looking at.

Thanks to $v=0$, the Hamiltonian equations reduce to:
$$
\left\{
\begin{array}{ll}
\ds \dot{g} & =\sum_{1}^{n-1}u_jY_j\\
\ds \dot{\lambda} & =\left(\sum_{1}^{n-1}u_j\ad(Y_j)\right)^*\lambda
\end{array}
\right.
$$

As the one-form $\lambda(t)$ is not zero but vanishes on $Y_1,\dots,Y_{n-1}$, that is on the left-invariant distribution $\Delta$, it can be written as $\lambda(t)=p(t)\omega$, where $p$ is a non-vanishing absolutely continuous function from $[0,T]$ into $\R$, and the ODE satisfied by $\lambda$ becomes
$$
\dot{p}(t)\omega=p(t)\left(\sum_{1}^{n-1}u_j\ad(Y_j)\right)^*\omega \ .
$$
This equality between one-forms is equivalent to the existence for almost every $t\in([0,T])$ of a real number $c(t)$ such that
\begin{enumerate}
\item[(a)] $\ds \qquad \left(\sum_{1}^{n-1}u_j\ad(Y_j)\right)^*\omega =c(t)\omega$
\item[(b)] $\qquad \dot{p}(t) =c(t)p(t)$
\end{enumerate}
Let us analyze the first equality. Let $Y_n\in\g$ such that $\omega(Y_n)=1$, and let $Y\in \Delta$.
Then
$$
\begin{array}{l}
\exists c\in\R \ \mbox{ s.t. } \ \ad(Y)^*\omega=c\omega \\
\qquad \Longleftrightarrow\  \exists c\in\R,\ \  \forall X\in\g \quad \langle\omega,\ad(Y)X\rangle=c\langle\omega,X\rangle\\
\qquad  \Longleftrightarrow\   \langle\omega,\ad(Y)Y_i\rangle=c\langle\omega,Y_i\rangle=0 \quad i=1,\dots,n-1\\
\qquad\qquad\qquad   \mbox{and }\ \langle\omega,\ad(Y)Y_n\rangle=c\langle\omega,Y_n\rangle=c\times 1=c\\
\qquad \Longleftrightarrow\  \ad(Y)\Delta\subseteq \Delta
\end{array}
$$

Let us set $\ds Y(t)=\sum u_i(t)Y_i$. As $g(t)$ is an abnormal curve, $Y(t)$ belongs to $\Delta$ and satisfies $\ad(Y(t))\Delta\subseteq \Delta$ for almost every $t\in [0,T]$. In other words $Y(t)$ satisfies:

\begin{equation}
\label{jesaispas} Y(t)\in \Delta\bigcap \mathcal{N}(\Delta)\quad \mbox{a.e. }t\in [0,T],
\end{equation}

where $\mathcal{N}(\Delta)$ stands for the normalizer of $\Delta$ in $\g$. Moreover
$$
c(t)=\langle\omega,\ad(Y(t))Y_n\rangle
$$
is  measurable essentially bounded as soon as $t\mapsto Y(t)$ is.

\vskip 0.2cm

Let us now look at the other condition, that is $g(t)$ belongs to $\Z$. As previously we denote by $\psi$ the function from $G$ to $\R$ defined by $\psi(g)=\langle\omega_g,\X_g\rangle=\langle\omega,F_g\rangle$. Recall that for all $Y\in\g$
$$
T_g\psi.TL_gY=\langle\omega,(D+\ad(F_g))Y\rangle \ .
$$
For $Y\in \Delta\bigcap \mathcal{N}(\Delta)$ and $g\in \Z$ we get
$$
T_g\psi.TL_gY=\langle\omega,DY\rangle \ .
$$
Indeed $\langle\omega,\ad(F(g))Y\rangle=-\langle\omega,\ad(Y)F(g)\rangle=0$ according to $F(g)\in\Delta$ and (\ref{jesaispas}). Consequently a curve $g(t)$ that verifies (\ref{jesaispas}) and $g(0)\in\Z$ is contained in $\Z$ if and only if $\langle\omega,D\dot{g}(t)\rangle=0$ almost everywhere, that is if and only if $\dot{g}(t)\in D^{-1}\Delta $ almost everywhere.

\vskip 0.2cm

Conversely let $g(t)$ be an absolutely continous curve defined on $[0,T]$ and such that $\dot{g}(t)=\dt g(t)$ belongs to
$$
\ab= \Delta\bigcap \mathcal{N}(\Delta)\bigcap D^{-1}\Delta
$$
for almost every $t$. It is straightforward to check that it is the projection of an abnormal extremal $(\lambda(t),g(t))$. The covector $\lambda(t)$ is equal to $p(t)\omega$, where $p(t)$ is the solution of the linear equation $\dot{p}=c(t)p$, with $p(0)\neq 0$ and $c(t)=\langle\omega,\ad(\dot{g}(t))Y_n\rangle$.
 
We are now in a position to state the following theorem.

\newtheorem{Anormales}[Zvariete]{Theorem}
\begin{Anormales} \label{Anormales}

The vector subspace $\ab$ of $\g$ defined by
$$
\ab=\Delta\bigcap \mathcal{N}(\Delta)\bigcap D^{-1}\Delta
$$
is a subalgebra of $\g$.

It generates a connected Lie subgroup of $G$ denoted by $A$, and for all $g\in \Z$ the coset $gA$ is included in the singular locus $\Z$.

The projections of the abnormal extremals are contained in $\Z$ and for $\g\in\Z$ the abnormal curves starting from $g$ are all the absolutely continuous curves contained in the coset $gA$.

Moreover the covecteur $\lambda(t)$ is up to a constant equal to $p(t)\omega$ where  $p(0)\neq 0$ and
$$
\dot{p}(t)=\langle\omega,\ad(\dot{g}(t))Y_n\rangle\ p(t)
$$
where $Y_n$ is such that $\langle\omega,Y_n\rangle=1$.

\end{Anormales}

\demo
We have only to prove that $\ab$ is a subalgebra, and that $gA$ is contained in $\Z$ as soon as $g\in\Z$.

An element $Y$ of $\g$ belongs to $\ab$ if and only if
$$
(i)\ \ Y\in\Delta \qquad\qquad  (ii)\ \ \ad(Y)\Delta\subset \Delta \qquad\qquad  (iii)\ \ DY\in\Delta.
$$

Let $Y,Z\in \ab$. By (ii) we get $[Y,Z]=\ad(Y)Z\in\Delta$.

Then $\ad(Y)\ad(Z)\Delta\subset \ad(Y)\Delta\subset \Delta$ according to (ii), hence $\ad([Y,Z])\Delta\subset \Delta$.

To finish $D[Y,Z]=[DY,Z]+[Y,DZ]\in \Delta$ since $DY$ (resp. $DZ$) belongs to $\Delta$ and $\ad(Z)$ ($\ad(Y)$) sends $\Delta$ into $\Delta$.

Consequently $[Y,Z]$ belongs to $\ab$ which is a subalgebra of $\g$.

To prove the second point let $g\in\Z$. Any absolutely continuous curve starting from $g$ and contained in $gA$, in particular any curve of the form $g\exp(tY)$ with $Y\in\ab$, satisfies the conditions discussed before the statement of the theorem. Consequently it is an abnormal curve, it is contained in $\Z$.

\hfill $\Box$

\vskip 0.3cm

\noindent{\bf Remarks}
\begin{enumerate}
\item
The dimension of $\ab$ is at most $n-2$. Indeed $\ab$ is included in $\Delta$ and its dimension cannot exceed $n-1$. Suppose it is equal to $n-1$, then $\Delta$ would be a subalgebra (because it would be included in $\mathcal{N}(\Delta)$), and would be moreover included in $D^{-1}\Delta$, hence invariant for $D$. But if $\Delta$ is a $D$-invariant subalgebra then the rank condition does not hold. Consequently $\dim(\ab)\leq n-2$.

On the other hand this dimension is equal to $n-2$ as soon as $\Delta$ is a subalgebra since in that case $\ab=\Delta\bigcap D^{-1}\Delta$, and $\Delta$ cannot be included in $D^{-1}\Delta$, which is therefore a codimension one subspace of $\g$ distinct from $\Delta$.
\item If $F(g)$ never belongs to the subalgebra $\ab$ (at no $g\in \Z$ or at least at no $g'\in gA$) then the controls associated to the abnormal extremals are uniquely defined and cannot involve the linear vector field. Consequently the abnormal curves are the left translation by $g$ of the geodesics of $A$ for the left-invariant metric induced by the one of the ARS.
\end{enumerate}


\section{Almost-Riemannian and sub-Riemannian structures}\label{ARSSRS}

In this section we show that a simple ARS on a $n$-dimensional Lie group is related to a codimension one regular sub-Riemannian structure on a $(n+1)$-dimensional Lie group. In other words the desingularization of a simple ARS is global and the desingularized structure is a classical sub-Riemannian one.

\subsection{The lift of a linear vector field}

Let $G$ be a connected Lie group, $\g$ its Lie algebra identified with the set of left-invariant vector fields, and $\X$ a linear vector field on $G$. The flow $(\varphi_t)_{t\in\R}$ of $\X$ is a one-parameter group of automorphisms of $G$, hence a Lie group morphism from $\R$ into $G$, so that we can consider the semi-direct product $\widetilde{G}=G\rtimes_{\varphi}\R$ defined by the law
$$
(g_1,\tau_1).(g_2,\tau_2)=(\varphi_{\tau_2}(g_1)g_2,\tau_1+\tau_2)
$$
The group $G$ can be identified with the closed and normal subgroup $G\times\{0\}$ of $\widetilde{G}$, and $\R$ with the closed subgroup $\{e\}\times \R$ of $\widetilde{G}$. This last is normal if and only if $\X=0$ (see \cite{AS01} for that construction).

\vskip 0.2cm

Let $\widetilde{\g}$ stand for the Lie algebra of $\widetilde{G}$. It is equal to some semi-direct product $\g\rtimes\R$ ; in view of its characterization the action of the Lie algebra $\R$ on $\g$ should be identified. Let $\widetilde{\X}$ denote the left-invariant vector field which is equal to $\ds (0,1)=\ddx{\tau}$ at the identity $\widetilde{e}=(e,0)$. This field generates the Lie algebra $\R$, so that $\exp(t\widetilde{\X})=(e,t)$. As it is left-invariant, we get also:
$$
\forall (g,\tau)\in\widetilde{G} \qquad (g,\tau)\exp(t\widetilde{\X})=(g,\tau).(e,t)=(\varphi_t(g),t+\tau),
$$
and
$$
\ds \widetilde{\X}(g,\tau)=\dto(\varphi_t(g),t+\tau)=(\X(g),1).
$$
Thus we get the equality $\ds \widetilde{\X}=\X+\ddx{\tau}$, and denoting as usual the derivation associated to $\X$ by $D$:
$$
\forall Y\in\g \qquad [\widetilde{\X},Y]=[\X,Y]=-DY.
$$
Finally the Lie algebra $\widetilde{\g}$ is the semi-direct product $\g\rtimes_{-D}\R$:
$$
[Y_1+\tau_1\widetilde{\X},Y_2+\tau_2\widetilde{\X}]=[Y_1,Y_2]-\tau_1DY_2+\tau_2DY_1.
$$


\subsection{Projection}

The group $G$ is diffeomorphic to the homogeneous space $\widetilde{G}/\R$, the set of right cosets of $\R$, but this last is not a quotient group in general because the subgroup $\R$ of $\widetilde{G}$ is not normal as soon as $\X\neq 0$. On the other hand the left-invariant vector fields of $\widetilde{G}$ can be projected to $\widetilde{G}/\R$ (this is always true for left-invariant vector fields and right cosets).

More accurately we get
$$
\begin{array}{l}
(e,\R)(g,\tau)\exp\left(t(Y,0)\right)=(e,\R)(g\exp(tY),\tau)=(g\exp(tY),\R)\\
(e,\R)(g,\tau)\exp(t\widetilde{\X})=(e,\R)(\varphi_t(g),\tau+t)=(\varphi_t(g),\R)
\end{array}
$$
so that we can identify the projection onto $\widetilde{G}/\R$ of the vector field $(Y,0)$ (resp. $\widetilde{\X}$) with the left-invariant vector field $Y$ (resp. the vector field $\X$) of $G$.


\subsection{Almost-Riemannian and sub-Riemannian structures}

Let us now consider an almost-Riemannian structure defined on $G$ by \\ $\{\X,Y_1,\dots,Y_{n-1}\}$, assumed to satisfy the rank condition, which implies
\begin{equation}
\mathcal{LA}\{\X,Y_1,\dots,Y_{n-1}\}=\R\X\oplus \g. \label{ARS34}
\end{equation}

On $\widetilde{G}$, we can define the left-invariant distribution
$$
\widetilde{\Delta}=\Span\{\widetilde{\X},Y_1,\dots,Y_{n-1}\}
$$
This distribution is left-invariant, of codimension 1, and on account of (\ref{ARS34}) we get
$$
\mathcal{LA}\{\widetilde{\X},Y_1,\dots,Y_{n-1}\}=\R\widetilde{\X}\oplus \g=\widetilde{\g}
$$
which shows that $\widetilde{\Delta}$ satisfies the rank condition. In order to define a corank one classical sub-Riemannian structure on $\widetilde{G}$ it remains to declare the vector fields $\widetilde{\X},Y_1,\dots,Y_{n-1}$ to be orthonormal.

Notice that the dynamics of the almost-Riemannian structure being $\dot{g}=v\X+u_1Y_1+\dots+v_{n-1}Y_{n-1}$ the one of the sub-Riemannian structure writes in the coordinates $(g,\tau)\in G\rtimes \R$:
$$
\left\{
\begin{array}{ll}
\dot{g} & =v\X+u_1Y_1+\dots+v_{n-1}Y_{n-1}\\
\dot{\tau} & =v
\end{array}
\right.
$$

\subsection{Optimality}

Let $u$ stand for the general control $(v,u_1,\dots,u_{n-1})$, let $g_0$ and $g_1$ be two points of $G$, and let $\bar{u}$ be a minimal control that steers $g_0$ to $g_1$ on the time interval $[0,T]$. By a minimal control is meant a control such that the associated curve $\gamma(t)$ minimizes the length between $g_0$ and $g_1$, and such that
$$
\norm{\bar{u}(t)}=\sqrt{\bar{v}^2(t)+\sum_{i=1}^{n-1}\bar{u}_i^2(t)}=\norm{\dot{\gamma}(t)} \qquad a.e.\ t\in [0,T].
$$
On $\widetilde{G}$ the control $\bar{u}$ steers $(g_0,\tau_0)$ to $(g_1,\tau_1)$ for any $\tau_0$ and $\tau_1$ that satisfy
\begin{equation}
\tau_1-\tau_0=\int_0^T \bar{v}(t)\ dt   \label{integrale}
\end{equation}
and it is minimal.

Indeed let $\tau\in\R$, let $\tilde{u}$ be a minimizing control that steers $(g_0,0)$ to $(g_1,\tau)$ in time $T$, and let $\widetilde{\gamma}$ be the associated trajectory. Its length is:
$$
l(\widetilde{\gamma})=\int_{0}^{T}\norm{\tilde{u}(t)}dt.
$$

Let $\gamma_1$ be the projection of $\widetilde{\gamma}$ on $G$. As it steers $g_0$ to $g_1$, we get:
$$
\int_{0}^{T}\norm{\bar{u}(t)}dt=l(\gamma)\leq l(\gamma_1)\leq l(\widetilde{\gamma})=
\int_{0}^{T}\norm{\tilde{u}(t)}dt.
$$ 

In conclusion the extremal trajectories on $G$ can be lifted to extremal trajectories on $\widetilde{G}$: the optimal control $\bar{u}$ belongs to the set of  optimal controls that steers $(g_0,0)$ to $(g_1,\tau)$ for some $\tau$, but this last is not known "à priori".


\section{Examples}\label{Examples}

\subsection{An example on the 2D affine group}

With the notations of Section \ref{ZAff} we consider the simple ARS defined by the left-invariant vector field $X$ and  the linear one $\X=Y+\mathcal{I}_*Y$, that is:
$$
\X(g)=\begin{pmatrix}x&y\\0&1\end{pmatrix}\begin{pmatrix}0&1\\0&0\end{pmatrix} - \begin{pmatrix}0&1\\0&0\end{pmatrix}\begin{pmatrix}x&y\\0&1\end{pmatrix}=\begin{pmatrix}0&x-1\\0&0\end{pmatrix}.
$$

The pair $\{\X,X\}$ is considered as an orthonormal frame, which defines the almost-Riemannian metric. Notice that
\begin{enumerate}
	\item[(i)] these two vector fields are linearly independant on the open and dense subset $\{x\neq 1\}$;
	\item[(i)] the rank condition is satisfied since $[X,\X]=Y$.
\end{enumerate}
Consequently the ARS is well defined and its singular locus is the line $\Z=\{x= 1\}$.

\vskip 0.2cm

The dynamic of that ARS is described by $\dot{g}=v\X+uX$, where $v,u\in\R$, or in coordinates:
$$
(\Sigma)\qquad\left\{\begin{array}{ll}
\dot{x} =& ux\\
\dot{y} =& v(x-1) 
\end{array} \right.
$$
The  associated Hamiltonian is
$$
\mathcal{H}=upx+vq(x-1)-\frac{1}{2}\nu(u^2+v^2), \qquad \nu\geq 0,
$$
where $(p,q)$ stands for the covecteur $\lambda$.


\subsubsection{Abnormal extremals}

We know by Theorem \ref{Anormales} that the abnormals are contained in some cosets of a subgroup of dimension at most $n-2=2-2=0$. Let us verify that on $Aff_+(2)$ the abnormals are fixed points. 

If $\nu=0$ we get
$$
\mathcal{H}_0=upx+vq(x-1)=u\left\langle\lambda,X \right\rangle +v\left\langle\lambda,\X\right\rangle
$$
so that $\mathcal{H}_0$ can have a maximum w.r.t. $(u,v)$ if and only if $\left\langle\lambda,X \right\rangle=\left\langle\lambda,\X\right\rangle=0$ (here $\left\langle .,. \right\rangle$ stands for the duality bracket). As $\lambda\neq 0$ the fields $X$ and $\X$ should be linearly dependant, which happens only in the singular locus. But then $\dot{y}=0$ and the abnormal extremals are fixed points.


\subsubsection{Normal extremals}
As usual, $\nu$ can be normalized to $\nu=1$. The computation of the partial derivatives of $\mathcal{H}$ w.r.t. $u$ and $v$ shows that $\mathcal{H}$ reaches its maximum for $u=px$ and $v=q(x-1)$, so that this maximum is equal to
$$
H=\frac{1}{2}p^2x^2+\frac{1}{2}q^2(x-1)^2,
$$
and the Hamiltonian equations are:
$$
\left\{\begin{array}{ll}
\dot{x} =& px^2\\
\dot{y} =& q(x-1)^2 
\end{array} \right.
\qquad
\left\{\begin{array}{ll}
\dot{p} =& -p^2x-q^2(x-1)\\
\dot{q} =& 0 
\end{array} \right.
$$

In the particular case where $q=0$ we get $p^2x^2=2H$, and setting $c=\sqrt{2H}$ we obtain $\dot{x}=\pm cx$ and $\dot{y}=0$.

Consequently the lines parametrized by $x(t)= x_0e^{\pm ct}$ and $y(t)=y_0$ are geodesics.


\subsubsection{Poincar\'e coordinates}
We are interested in the geodesics starting from  the singular locus, that is with initial conditions $(x_0=1,y_0)$.

We can restrict ourselves to the geodesics parametrized by arclength, that is to $\ds H=\frac{1}{2}$. As $\ds H=\frac{1}{2}p_0^2$ at $t=0$ for geodesics starting from $(1,y_0)$, this amounts to take $p_0=\epsilon=\pm 1$. Moreover we can also assume $q=q_0>0$ because of the symmetry w.r.t. the $x$ axis.

To solve the Hamiltonian equations two changes of variables are done. Since $H$ is constant we can first set:
$$
\left\{\begin{array}{ll}
\epsilon\cos(\alpha) & = px\\
\epsilon\sin(\alpha) & = q(x-1)
\end{array} \right.
$$
Derivating the first of these formulas we get
$$
-\epsilon\dot{\alpha}\sin(\alpha)=p\dot{x}+\dot{p}x=p^2x^2-p^2x^2-q^2x(x-1)=-qx\epsilon\sin(\alpha),
$$
so that
$$
\dot{\alpha}=qx=\epsilon\sin(\alpha)+q.
$$
For an initial point in the singular locus, that is $x_0=1$, and according to $q\neq 0$, we get $\alpha_0=0$ (or $k\pi$). The equation to be solved is consequently
\begin{equation}\label{edoalpha}
\dot{\alpha}=p_0\sin(\alpha)+q \quad \mbox{ with }\ \alpha(0)=0 \ \mbox{ and }\ q> 0.
\end{equation}

In order to solve (\ref{edoalpha}) we do the change of variable $\ds\tau=\tan(\frac{\alpha}{2})$, and the equation becomes:
\begin{equation}\label{edotau}
\ds \dot{\tau}=\epsilon\tau+\frac{1}{2}q(1+\tau^2)\qquad \tau(0)=0.
\end{equation}
Before solving (\ref{edotau}), we can notice that a straightforward computation gives:
$$
\left\{
\begin{array}{l}
\ds x=1+\frac{\epsilon}{q}\frac{2\tau}{1+\tau^2}\\
\\
\ds y=qt+\frac{\epsilon}{q}\frac{2\tau^2}{1+\tau^2}-2\arctan(\tau)+y_0
\end{array}
\right.
$$
The separation of the variables in $(\ref{edotau})$ leads to
$$
t=\frac{2}{q}\int \frac{d\tau}{\tau^2+2r\tau+1}
$$
where $\ds r=\frac{\epsilon}{q}$. To compute these primitives we should distinguish three cases according to the sign of the discriminant $4(r^2-1)$ of $\tau^2+2r\tau+1$.

\vskip 0.2cm

\noindent \textbf{Case 1} ($0<q<1$)

We get
$$
\ds \tau(t)=\frac{q\left(\exp(\sqrt{1-q^2}\ t )-1\right)}{\epsilon+\sqrt{1-q^2}+\left(\sqrt{1-q^2}-\epsilon\right)\exp(\sqrt{1-q^2}\ t)}.
$$
We are interested in the first return to the singular locus $\{x=1\}$. It happens when $\tau(t)$ vanishes for some $t_0$, or as a limit if $\ds \tau(t)\longmapsto_{t\mapsto t_0}\pm \infty$.

For $0<q<1$ the function $\tau(t)$ does not vanish for $t>0$ but tends to $+\infty$ when $t$ tends to $\ds \frac{1}{\sqrt{1-q^2}}\ln\left(\frac{\epsilon+\sqrt{1-q^2}}{\epsilon-\sqrt{1-q^2}}  \right)$, which is positive only for $\epsilon=1$. Consequently when $\epsilon=1$ and $t$ tends to $\ds \frac{1}{\sqrt{1-q^2}}\ln\left(\frac{1+\sqrt{1-q^2}}{1-\sqrt{1-q^2}}  \right)$, then $x(t)$ tends to $1$, and
$$
y(t)\longmapsto \frac{q}{\sqrt{1-q^2}}\ln\left(\frac{1+\sqrt{1-q^2}}{1-\sqrt{1-q^2}}  \right)+\frac{2}{q}-\pi+y_0.
$$
Notice that the geodesics do not go back to the $y$-axis in positive time for $\epsilon=-1$, that is for $x\leq 1$, as illustrated on the following picture:

\includegraphics{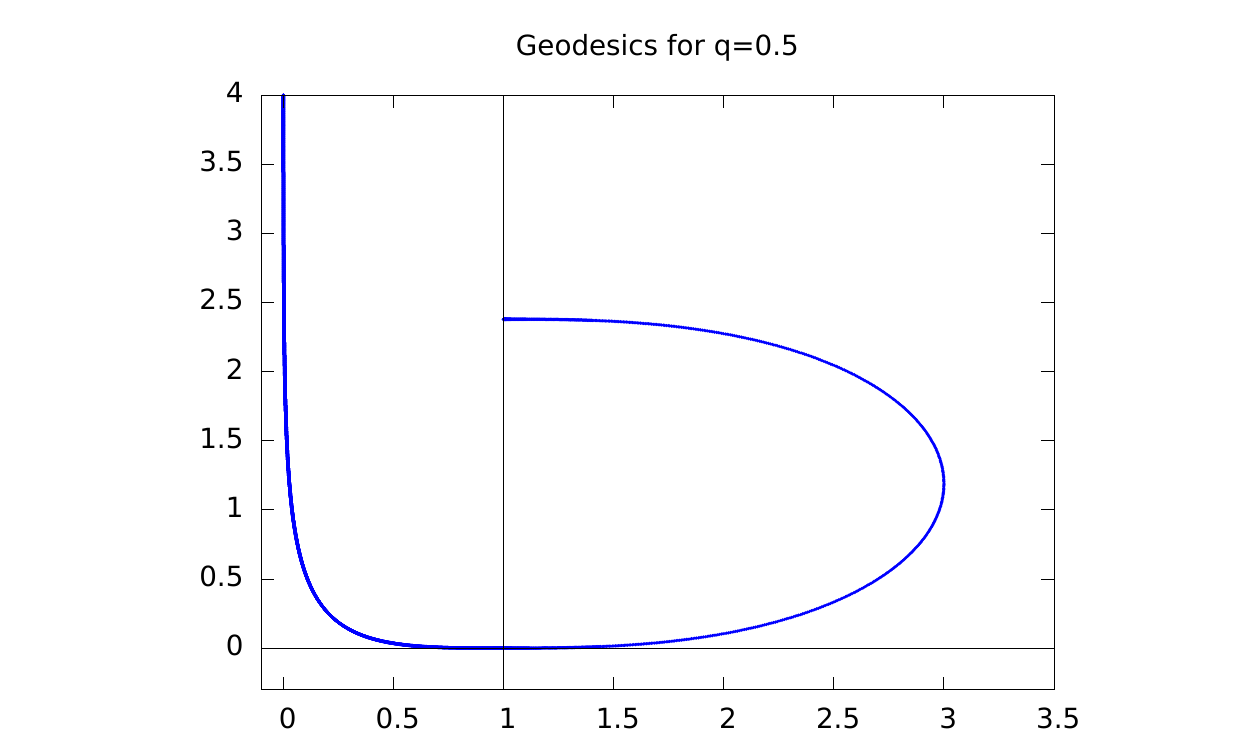}

\vskip 0.2cm

\noindent \textbf{Case 2} ($q=1$)

We get $\ds \tau(t)=\frac{t}{2-\epsilon t}$ hence
$$
\left\{
\begin{array}{l}
\ds x(t)=1+\frac{2\epsilon t-t^2}{2-2\epsilon t+t^2}\\
\\
\ds y(t)=t+\frac{\epsilon t^2}{2-2\epsilon t+t^2}-2\arctan(\frac{t}{2-\epsilon t})+y_0
\end{array}
\right.
$$
The trajectory belongs to the singular locus, that is  $x(t)=1$, for $t=0$ and $t=2\epsilon$.

When $t\mapsto 2\epsilon$, $y(t)$ tends to $\epsilon(4-\pi)$, that is the locus of the first return to the $y$-axis. Notice that as well as in the first case the geodesics do not go back to the $y$-axis in positive time for $\epsilon=-1$, that is for $x\leq 1$.

\includegraphics{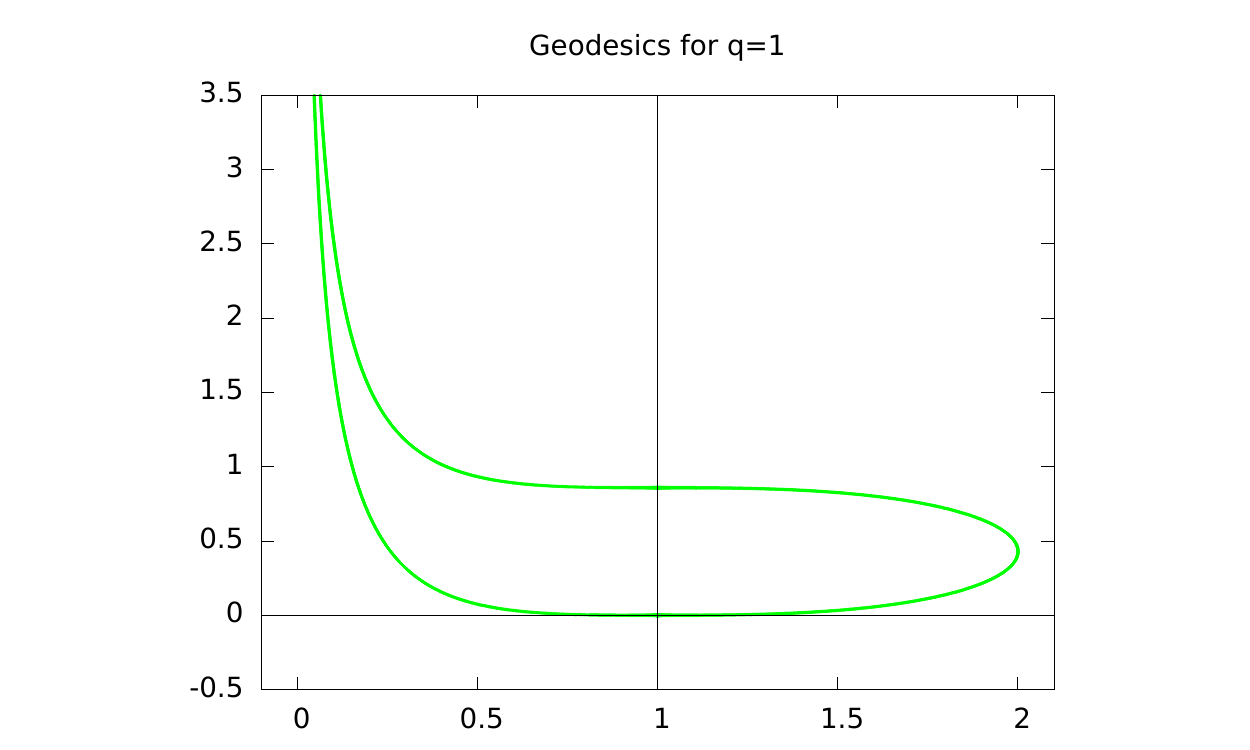}

\vskip 0.2cm

\noindent \textbf{Case 3} ($q>1$)

We get
$$
\ds \tau(t)=\frac{\sqrt{q^2-1}}{q}\tan\left(\frac{1}{2}t\sqrt{q^2-1}+\epsilon\theta\right)-\frac{\epsilon}{q}
$$
where $\ds \theta=\arctan\left(\frac{1}{\sqrt{q^2-1}}\right)$.
 
As in the first case $\tau(t)$ does not vanish, hence $x(t)$ cannot be equal to 1, for $t>0$, but  $\ds \tau(t)\mapsto +\infty$ when
$
\ds t\mapsto\frac{\pi-2\epsilon \theta}{\sqrt{q^2-1}}
$.
Consequently when $t$ tends to that limit we get:
$$
\left\{
\begin{array}{ll}
\ds x(t)=1+\frac{\epsilon}{q}\frac{2\tau}{1+\tau^2} & \longmapsto 1\\
&\\
\ds y(t)=qt+\frac{\epsilon}{q}\frac{2\tau^2}{1+\tau^2}-2\arctan(\tau) +y_0& \ds \longmapsto q\frac{\pi-2\epsilon \theta}{\sqrt{q^2-1}}+\frac{2\epsilon}{q}-\pi+y_0 
\end{array}
\right.
$$
The behaviour is here different from the two previous cases because the geodesics return to the axis $\{x=1\}$ in both cases $\epsilon=\pm 1$.

\includegraphics{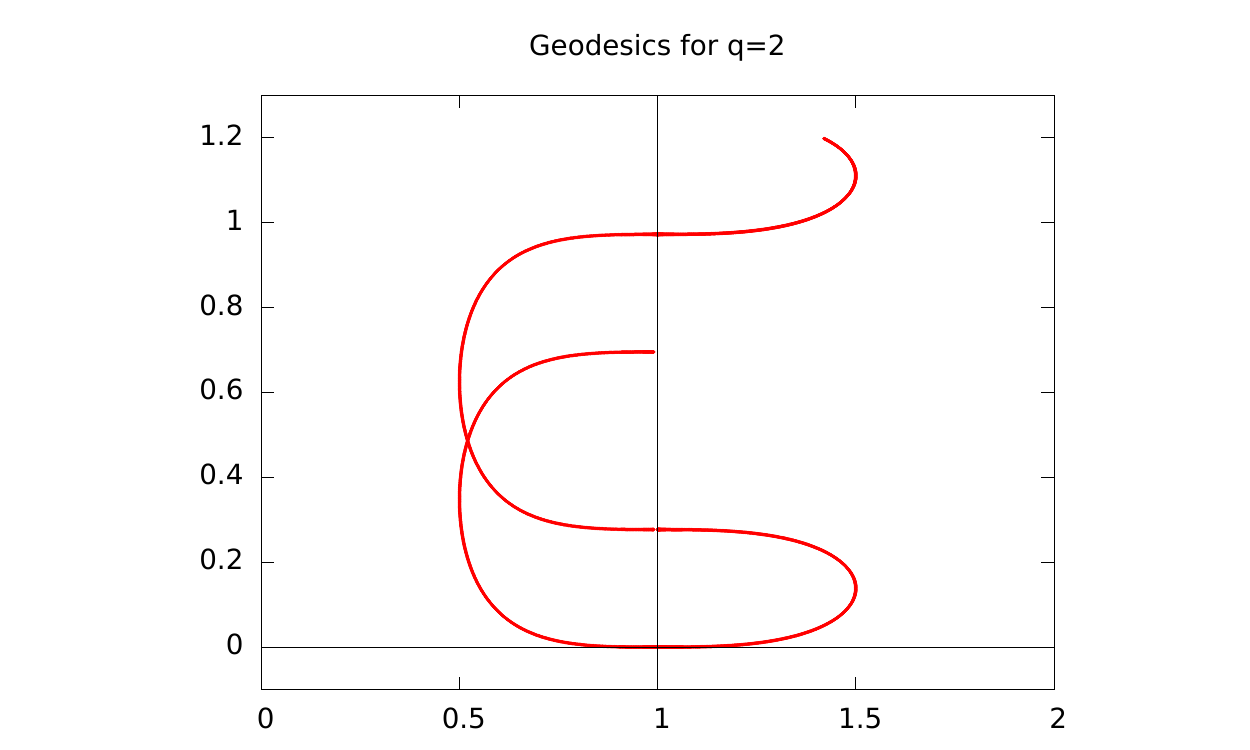}

To finish the three behaviours are compared in a unique picture and the unit sphere is drawn

\includegraphics[width=0.5\linewidth]{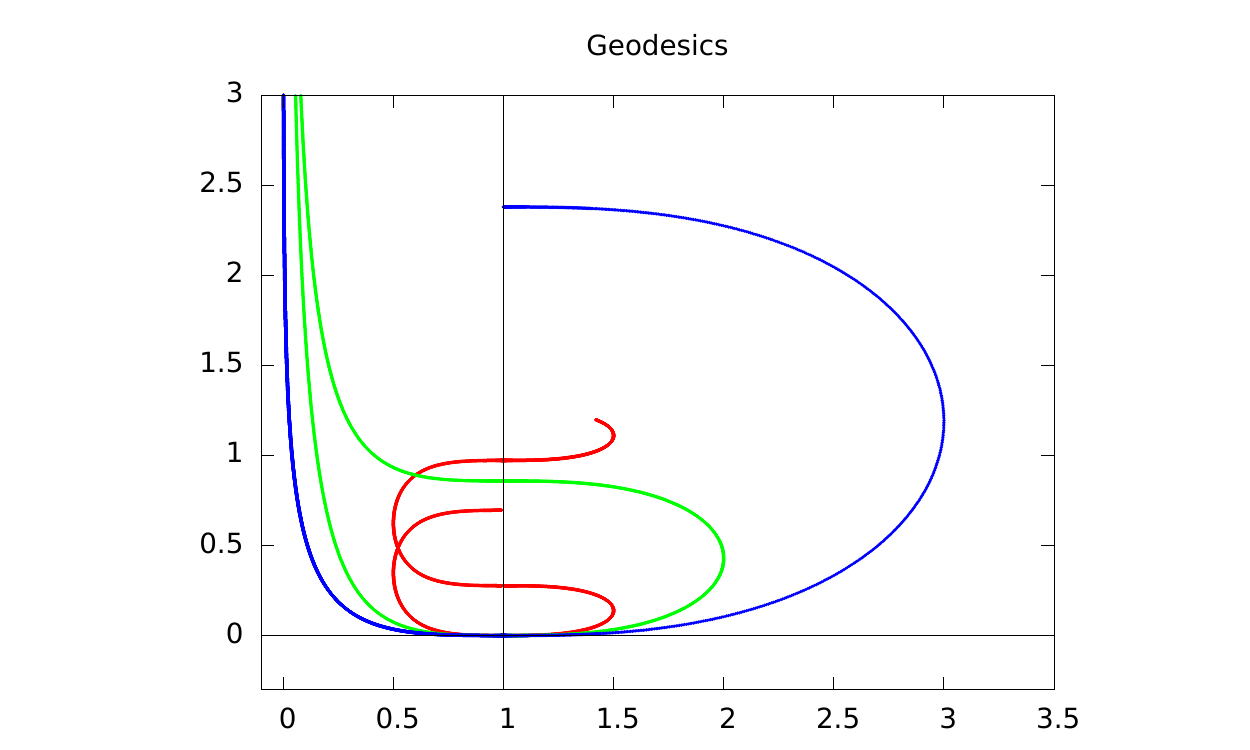}
\includegraphics[width=0.5\linewidth]{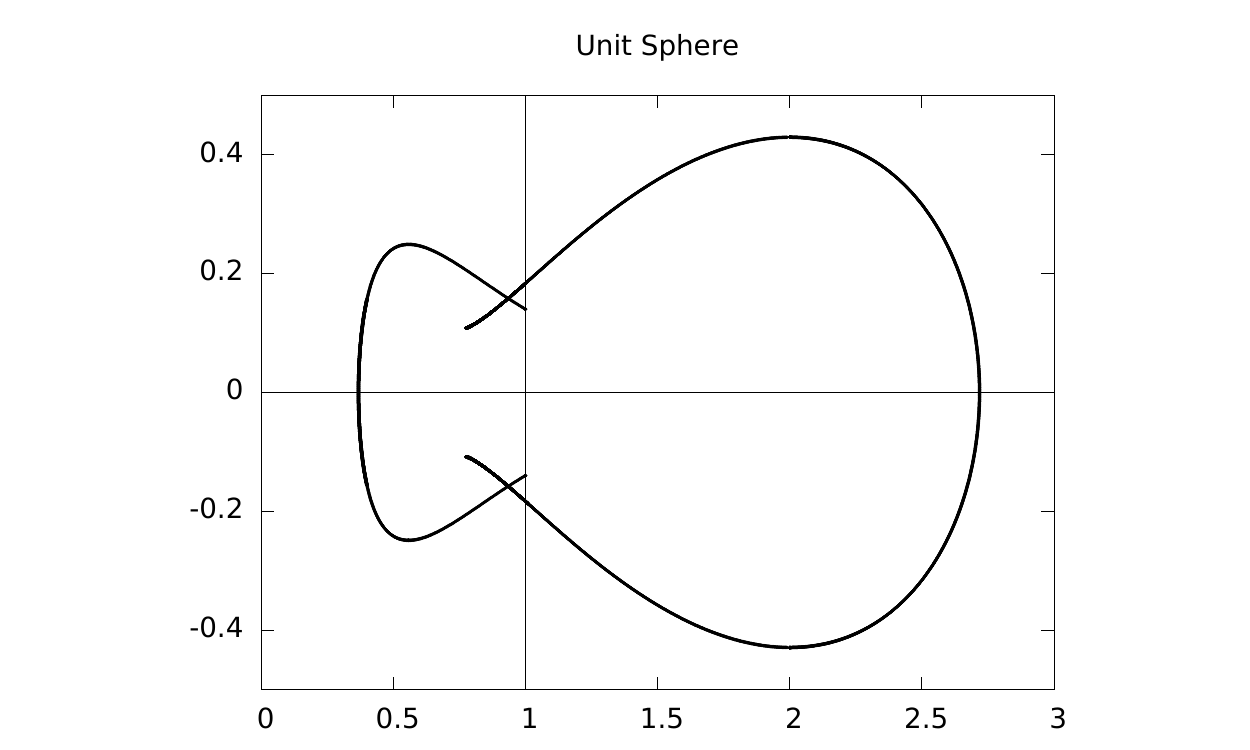}


\subsection{An example on the Heisenberg group}
The equations of the following example are only stated, not solved. Actually the purpose is here to illustrate some of the theoretical results, in particular the desingularization process.

With the notations of Section \ref{ZHeisenberg} we consider on the Heisenberg group $G$ the ARS defined by the set of vector fields $\{\X,X,Z\}$ where
$$
\ds \X=x\ddx{y}+\frac{1}{2}x^2\ddx{z} \quad
\mbox{is associated to the derivation } D=\begin{pmatrix}0&0&0\\1&0&0\\0&0&0\end{pmatrix}
$$

Here $\Delta=\mbox{Span}\{X,Z\}$ is an ideal of $\g$, and according to Proposition \ref{Omegaclosed} the singular locus is a subgroup. An easy computation shows that it is equal to the set of fixed points of $\X$: $\Z=\{x=0\}$. On the other hand the rank condition is satisfied since $[X,\X]=DX=Y$, and the ARS is well defined.

The dynamic of that ARS is here $\dot{g}=v\X+u_1Y_1+u_2Y_2$, where $v,u_1,u_2\in\R$, or in coordinates
$$
(\Sigma)\qquad\left\{\begin{array}{ll}
\dot{x} & =u_1\\
\dot{y} & =vx \\
\dot{z} & =u_2+\frac{1}{2}vx^2
\end{array} \right.
$$
and the associated Hamiltonian is
$$
\begin{array}{ll}
\mathcal{H} &=\left\langle \lambda,v\X+u_1Y_1+u_2Y_2\right\rangle-\frac{1}{2}\nu(v^2+u_1^2+u_2^2)\\
 & =pu_1+qvx+r(u_2+\frac{1}{2}vx^2)-\frac{1}{2}\nu(v^2+u_1^2+u_2^2), \qquad \nu\geq 0.
\end{array}
$$
where $\lambda=(p,q,r)$.


\subsubsection{Abnormal extremals}

Let us compute the abnormals extremals. For $\nu=0$ we get
$$
\mathcal{H}_0=u_1p+u_2r+v(qx+\frac{1}{2}rx^2)
$$
so that $\mathcal{H}_0$ presents a maximum if and only if $p=r=qx=0$. But $p=r=0$ imply $q\neq 0$, hence $x=0$. Moreover $q$ is a constant because $-\ddx{y}\mathcal{H}_0=0$, and choosing $q=1$ we finally obtain $x=0$ (that is the abnormals are contained in the singular locus) together with the equations
$$
\left\{\begin{array}{ll}
\dot{y} & =0 \\
\dot{z} & =u_2
\end{array}\right.
$$
In other words the vertical lines through a point $(0,y_0,0)$ are abnormal extremals.

On the other hand the normalizer of $\Delta=\Span\{X,Z\}$ being $\g$ and $D^{-1}\Delta$ being equal to $\Span\{Y,Z\}$, we know by Theorem \ref{Anormales} that the abnormals extremals are the curves contained in $\Z=\{x=0\}$ and parallel to $\ab=\Span\{Y\}$, that is the computed result.


\subsubsection{Normal extremals}
As previously $\nu$ is normalized to $\nu=1$. The computation of the partial derivatives of $\mathcal{H}_1$ w.r.t. $v$, $u_1$ and $u_2$ shows that $\mathcal{H}_1$ reaches its maximum for $u_1=p$, $u_2=r$, and $v=qx+\frac{1}{2}rx^2$, so that this maximum is equal to
$$
\ds H=\frac{1}{2}\left(p^2+r^2+(qx+\frac{1}{2}rx^2)^2\right),
$$
and the Hamiltonian equations are:
$$
\left\{\begin{array}{ll}
\dot{x} & = p\\
\dot{y} & = qx^2+\frac{1}{2}rx^3\\
\dot{z} & = r+\frac{1}{2}qx^3+\frac{1}{4}rx^4
\end{array} \right.
\qquad
\left\{\begin{array}{ll}
\dot{p} & = -(q+rx)(qx+\frac{1}{2}rx^2)\\
\dot{q} & = 0 \\
\dot{r} & =0
\end{array} \right.
$$

One can see at once that the abnormal extremals are normal, they are obtained for $x=p=0$.

On the other hand the case $q=r=0$, $p\neq 0$ gives $\dot{x}=p$ and $\dot{p}=\dot{y}=\dot{z}=0$, hence all the lines parallel to the $x$-axis.


\subsubsection{Poincar\'e coordinates}

Apart from these simple cases, we are lead to isolate the two following equations:
$$
\left\{\begin{array}{ll}
\dot{x} & = p\\
\dot{p} & = -(q+rx)(qx+\frac{1}{2}rx^2)
\end{array} \right.
$$

We know that $\ds H=\frac{1}{2}\left(p^2+r^2+qx+(qx+\frac{1}{2}rx^2)^2\right)$ is constant along the extremals and, as $r$ is also constant, we can set
$$
p^2+(qx+\frac{1}{2}rx^2)^2=2H^2-r^2=c^2,
$$
and then
$$
\left\{\begin{array}{ll}
c\cos\alpha =& p\\
c\sin\alpha =& qx+\frac{1}{2}rx^2
\end{array} \right.
$$

By a first derivation we get
$$
\left\{\begin{array}{ll}
-c\dot{\alpha}\sin\alpha & = \dot{p}=-(q+rx)(qx+\frac{1}{2}rx^2)=-(q+rx)c\sin\alpha\\
\ \ c\dot{\alpha}\cos\alpha & = (q+rx)\dot{x}=(q+rx)p=(q+rx)c\cos\alpha
\end{array} \right.
$$
so that $\dot{\alpha}=q+rx$. A second derivation gives:
$$
\ddot{\alpha}=r\dot{x}=rp=cr\cos\alpha,
$$
that is the pendulum equation.

For an initial point in the singular locus, that is $x_0=0$, we get $\alpha_0=0$ (or $k\pi$), and $p_0=\pm c$. The equation to be solved is consequently
$$
\ddot{\alpha}=p_0r\cos(\alpha) \quad \mbox{ with $\alpha(0)=0$ .}
$$


\subsubsection{Desingularization}

The Lie algebra generated by $\X,Y_1,Y_2$ is equal to $Sp\{\X,X,Y,Z\}=\g\oplus\R\X$. This algebra is the Engel one since
$[X,\X]=Y$, $[X,Y]=Z$, and the other brackets vanish. According to Section \ref{ARSSRS} it it isomorphic to the Lie algebra of the semidirect product of $G$ by $\R$, that is $\g\oplus\R\widetilde{\X}$ where $\ds \widetilde{\X}=\X+\ddx{w}$. The lift of the almost-Riemannian distribution to $\widetilde{G}=G\rtimes \R$ is the dimension 3 classical sub-Riemannian distribution $\Delta=\Span\{\widetilde{\X},Y_1,Y_2\}$, the dynamics of which writes in natural coordinates:
$$
\left\{\begin{array}{ll}
\dot{x} & =u_1\\
\dot{y} & =vx \\
\dot{z} & =u_2+\frac{1}{2}vx^2\\
\dot{w} & =v
\end{array} \right.
$$
The associated Hamiltonian is
$$
\ds \widetilde{\mathcal{H}}=\mathcal{H}+sv=u_1p+ u_2r+v(s+qx+\frac{1}{2}rx^2)-\frac{1}{2}\nu(v^2+u_1^2+u_2^2), \qquad \nu\geq 0.
$$
where $\widetilde{\lambda}=(p,q,r,s)\in \widetilde{\g}^*$.

\vskip 0.2cm

\noindent \textit{Abnormal extremals}. An easy computation shows first that $\nu=0$ implies $p=r=s+qr=0$, and then that the abnormal trajectories are the lines parallel to $Oz$ (after a distinction between the cases $s=0$ and $s\neq 0$).

\vskip 0.2cm

\noindent \textit{Normal extremals}. Setting $\nu=0$ we obtain at once the maximized Hamiltonian
$$
\ds \widetilde{H}=\frac{1}{2}\left(p^2+r^2+(s+qx+\frac{1}{2}rx^2)^2\right),
$$
and the Hamiltonian equations:
$$
\left\{\begin{array}{ll}
\dot{x} & = p\\
\dot{y} & = sx+qx^2+\frac{1}{2}rx^3\\
\dot{z} & = r+\frac{1}{2}sx^2+\frac{1}{2}qx^3+\frac{1}{4}rx^4\\
\dot{w} & = s+qx+\frac{1}{2}rx^2
\end{array} \right.
\qquad
\left\{\begin{array}{l}
\dot{p}= -(q+rx)(s+qx+\frac{1}{2}rx^2)\\
\dot{q} = \dot{r}=\dot{s}=0 
\end{array} \right.
$$
The almost-Riemannian dynamics is the part of that last obtained for $s=0$.


\section{Extension of the model and Equivalence}\label{EME}

The model analyzed in the previous sections can be generalized in two ways. Firstly on a connected  $n$-dimensional Lie group $G$, instead of considering $n-1$ left-invariant vector fields and one linear, we can consider a set $\{F_1,\dots,F_n\}$ of $n$ affine vector fields. By affine vector field is meant an element of the normalizer of the Lie algebra $\g$ of $G$ in the Lie algebra of analytic vector fields on $G$. The affine vector fields are sums of right-invariant, left-invariant and linear vector fields (see \cite{Jou09}), so that they are either left-invariant or the sum of a linear vector field and a right-invariant one. 

One of the reasons for considering affine vector fields instead of linear ones is that it may happen that the zero locus of several non left-invariant vector fields contain no common point, which would be the case, the common point being the identity $e$, if all these fields were linear.

\newtheorem{generaldef}[DefARS]{Definition}
\begin{generaldef} \label{generaldef}
A \textbf{general almost-Riemannian structure} on a connected $n$-dimensional Lie group $G$ is defined by a set $\{F_1,\dots,F_n\}$ of $n$ affine vector fields that satisfy the two following properties:
\begin{enumerate}
	\item[(i)] The rank of the distribution defined by $\Span\{F_1,\dots,F_n\}$ is equal to $n$ on a proper open subset $U$ of $G$.
	\item[(ii)] This distribution satisfies the rank condition.
\end{enumerate}
The almost-Riemannian metric is defined by declaring the set $\{F_1,\dots,F_n\}$ of vector fields to be an orthonormal frame.
\end{generaldef}

The second generalization consists in considering the same framework on homogeneous spaces. Let $H$ be a closed subgroup of $G$, and $G/H$ the manifold of right cosets of $H$. It is a well known fact that the left-invariant vector fields can be projected to the quotient group. The set of their projections will be denoted by $\Pi_*\g$; it is a Lie algebra. On the other hand an affine vector field on $G/H$ is defined in a natural way as the projection, whenever it exists, of an affine vector field of $G$. However it is shown in \cite{Jou09} that the affine vector fields of $G/H$ are exactly the elements of the normalizer of $\Pi_*\g$ in the Lie algebra of analytic vector fields on $G/H$. Consequently we can copy the previous definition:

\newtheorem{homdef}[DefARS]{Definition}
\begin{homdef} \label{homdef}
A \textbf{general almost-Riemannian structure} on a connected $n$-dimensional homogeneous space $G/H$ is defined by a set $\{F_1,\dots,F_n\}$ of $n$ affine vector fields that satisfy the two following properties:
\begin{enumerate}
	\item[(i)] The rank of the distribution defined by $\Span\{F_1,\dots,F_n\}$ is equal to $n$ on a proper open subset $U$ of $G/H$.
	\item[(ii)] This distribution satisfies the rank condition.
\end{enumerate}
The almost-Riemannian metric is defined by declaring the set $\{F_1,\dots,F_n\}$ of vector fields to be an orthonormal frame.
\end{homdef}

\noindent \textbf{Remark}.
In both cases all the vector fields are analytic and the condition that the rank of the distribution defined by $\Span\{F_1,\dots,F_n\}$ be equal to $n$ on an proper open subset $U$ implies that $U$ is dense in $G$ (or $G/H$).
	
\vskip 0.2cm

Consider now a connected and $n$-dimensional manifold $M$ endowed with an almost-Riemannian structure defined by a set $\{f_1,\dots,f_n\}$ of vector fields that verifies
\begin{enumerate}
	\item[(i)] The rank of $\Span\{f_1,\dots,f_n\}$ is equal to $n$ on an open and dense subset $U\subsetneq M$.
	\item[(ii)] The set $\{f_1,\dots,f_n\}$ satisfies the rank condition.
\end{enumerate}
As previously the almost-Riemannian metric is defined by considering the vector fields $\{f_1,\dots,f_n\}$ as an orthonormal frame.

If we want this structure to be equivalent to a general almost-Riemannian structure on a Lie group or a homogeneous space, we should add the two following necessary conditions:
\begin{enumerate}
	\item[(iii)] The vector fields $f_1,\dots,f_n$ generate a finite dimensional Lie algebra.
	\item[(iv)] The  vector fields $f_1,\dots,f_n$ are complete.
\end{enumerate}

These conditions turn out to be sufficient.

\newtheorem{ARSequivalence}[Zvariete]{Theorem}
\begin{ARSequivalence} \label{ARSequivalence}
Under the conditions $(i)$ to $(iv)$ the manifold $M$ is diffeomorphic to a homogeneous space $G/H$ and the vector fields $f_1,\dots,f_n$ are related by this diffeomorphism to invariant or affine vector fields on $G/H$ so that the almost-Riemannian structure on $M$ is equivalent to a general almost-Riemannian structure on $G/H$.

If condition $(iv)$ is removed, then the statement holds locally.
\end{ARSequivalence}

\demo 
Let $\Li$ be the (finite dimensional) Lie algebra generated by $\{f_1,\dots,f_n\}$ and for each subset $I$ of $\{1,2,\dots,n\}$ let us denote by $\Li(I)$ the ideal  generated in $\Li$ by $\{f_i;\ i\in I\}$. Let us choose $I_0$ as one of the subsets of $\{1,2,\dots,n\}$ that satisfy the rank condition and whose cardinal is minimal among the subsets of $\{1,2,\dots,n\}$ that satisfy that condition. This is possible since on the one hand $\Li=\Li(\{1,2,\dots,n\})$ satisfies the rank condition, and on the other hand the choice is done in a finite set.

We set $\Lo=\Li(I_0)$. Let $G$ be the simply connected Lie group whose Lie algebra $\g$ is isomorphic to $\Lo$. According to the Equivalence Theorem of \cite{Jou09} there exist a closed subgroup $H$ of $G$ and a diffeomorphism $\Phi$ from $G/H$ onto $M$ such that
$$
\Lo=\Phi_*(\Pi_*\g),
$$
where $\Pi_*\g$ stands for the Lie algebra of invariant vector fields of $G/H$.

We can assume without loss of generality that $I_0=\{1,\dots,k\}$ with $k<n$. As $\Lo$ is an ideal of $\Li$, the fields $f_{k+1},\dots,f_n$ belong to the normalizer of $\Lo$. Their images $F_i=\Phi^{-1}_*f_i$, $i=k+1,\dots,n$ belong to the normalizer of $\Pi_*\g$, hence are affine.

If we do not assume the $f_i$'s to be complete then the Equivalence Theorem applies locally. This is not explicitly stated in \cite{Jou09}, but is clear from the proof of the Equivalence Theorem (Theorem 5.1 of \cite{Jou09}).

\hfill $\Box$

\newtheorem{ARSLieequi}[localsubgroup]{Corollary}
\begin{ARSLieequi} \label{ARSLieequi}
Under the previous conditions, if moreover $M$ is simply connected and $\dim(\Lo)=\dim(M)$, then the ARS on $M$ is diffeomorphic to a general ARS on a Lie group.

If $M$ is not simply connected, or if the vector fields $f_i$ are not complete then the same statement holds locally.
\end{ARSLieequi}

\demo
With the notations of the previous proof we get: if $\dim(\Lo)=\dim(M)$ then $\dim(M)=\dim(G)$, the subgroup $H$ is discrete, hence reduced to the identity because of the simply connectedness assumption on $M$.

\hfill $\Box$

\vskip 0.2cm

\noindent \textbf{Remark} (Due to Jean-Paul-Gauthier)

In the statement of Theorem \ref{ARSequivalence} the ARS on the manifold $M$ is globally defined by $n$ orthonormal vector fields that generate a finite dimensional Lie algebra. The problem to know if a general ARS can be defined, at least locally, by a finite set of orthonormal vector fields that generate a finite dimensional Lie algebra is open.


\section{Appendix}

Let $\X$ be a linear vector field on a connected Lie group $G$. We denote by $F$ the mapping from $G$ into $\g$ defined by
\begin{equation}
\label{Xitranslated} F(g)=TL_{g^{-1}}\X_g
\end{equation}

\newtheorem{A1}[Omegaclosed]{Proposition}
\begin{A1} \label{A1}
\begin{enumerate}
\item For all $Y\in\g$
\begin{equation}\label{TeF}
\ds \dto F(\exp(tY))=DY
\end{equation}
\item 
For all $g\in G$ and for all $Y\in\g$
\begin{equation}\label{FgExp}
F(g\exp tY)=F(\exp tY)+e^{-t\ad(Y)}F(g)
\end{equation}
and
\begin{equation}\label{TgF}
T_gF=(D+\ad(F_g))\circ TL_{g^{-1}}
\end{equation}
\end{enumerate}
\end{A1}

\demo
\begin{enumerate}
\item Denoting by $(\varphi_t)_{t\in\R}$ the flow of $\X$ and according to Formula (\ref{equationfonda}) (see Section \ref{BD}) we have:
$$
\begin{array}{ll}
\ds\dto F(\exp(tY)) &\ds =\dto TL_{(\exp(tY))^{-1}}\X_{\exp(tY)}\\
\ds                  &\ds =\dto \dso \exp(-tY)\varphi_s(\exp(tY))\\
\ds                  &\ds =\dso \dto \exp(-tY)\exp(te^{sD}Y) \\
\ds                  &\ds =\dso(-Y+e^{sD}Y)=DY.
\end{array}
$$
\item Let $Y\in\g$. According to Formula \ref{Bourbak} of Section \ref{BD} we get:
$$
\begin{array}{ll}
\ds F(g\exp(tY)) & =\ds TL_{(g\exp(tY))^{-1}}\X_{g\exp(tY)}\\
\ds                  &\ds = TL_{\exp(-tY)}TL_{g-1}(TL_g\X_{\exp(tY)}+TR_{\exp(tY)}\X_g)\\
\ds                  &\ds =F(\exp(tY))+\Ad(\exp(-tY))F(g)\\
\ds                  &\ds =F(\exp(tY))+e^{-t\ad(Y)}F(g),
\end{array}
$$
that is Formula (\ref{FgExp}). Derivating that formula at $t=0$ we obtain at once Formula (\ref{TgF}):
$$
\begin{array}{ll}
T_gF.Y_g & \ds =T_gF. TL_g.Y=\dto F(g\exp(tY))\\
         & \ds =\dto\left(F(\exp(tY))+e^{-t\ad(Y)}F(g)\right)\\
         & \ds =DY-\ad(Y)F(g)=DY+\ad(F(g))Y.
\end{array}
$$
\end{enumerate}
\hfill $\Box$

The next very useful corollary is immediate from Formula (\ref{TgF}).
\newtheorem{C1}[localsubgroup]{Corollary}
\begin{C1} \label{C1}
Let $\omega$ be a left-invariant one-form on $G$, and $\psi$ the real valued function defined on $G$ by $\psi(g)=\langle\omega,F(g)\rangle$. Then
$$
\forall Y\in \g \quad T_g\psi.Y_g=\langle\omega,DY+\ad(F_g)Y\rangle.
$$
This can also be written $\ds T_g\psi\circ TL_g=\left( DY+\ad(F_g)\right)^*\omega$.
\end{C1}

To finish we compute $F$ as a power series.

\newtheorem{A2}[Omegaclosed]{Proposition}
\begin{A2} \label{A2}
For all $g\in G$, $Y\in\g$, and $t\in\R$:
\begin{enumerate}
\item[(i)]
$
\ds \forall k\geq 0 \qquad \frac{d^k}{dt^k}F(g\exp tY)=(-1)^{k-1}\ad^{k-1}(Y)DY+(-1)^k\ad^k(Y)F(g\exp tY)
$
\item[(ii)] $\ds F(\exp tY)=\sum_{k=1}^{+\infty}(-1)^{k-1}\frac{t^k}{k!}\ad^{k-1}(Y)DY$
\end{enumerate}
\end{A2}

\demo
\begin{enumerate}
\item[(i)] According to Formula (\ref{TgF}) of the previous proposition we get
$$
\dt F(g\exp(tY))=T_{g\exp(tY)}F.Y_{g\exp(tY)}=DY-\ad(Y)F(g\exp(tY)),
$$
that is the expected formula for $k=1$. By induction we have also:
$$
\begin{array}{ll}
\ds\frac{d^{k+1}}{dt^{k+1}}F(g\exp tY) & =(-1)^k\ad^k(Y)(DY-\ad(Y)F(g\exp tY))\\
                                       & =(-1)^k\ad^k(Y)DY+(-1)^{k+1}\ad^{k+1}(Y)F(g\exp tY).
\end{array}
$$
\item[(ii)] We have $F(\exp(0Y))=F(e)=0$ and for $k\geq 1$ the formulas of item $(i)$ applied to $g=e$ and $t=0$ gives
$$
 \frac{d^k}{dt^k}_{|t=0}F(\exp(tY))=(-1)^{k-1}\ad^{k-1}(Y)DY,
$$
so that by analycity
$\ds F(\exp tY)=\sum_{k=1}^{+\infty}(-1)^{k-1}\frac{t^k}{k!}\ad^{k-1}(Y)DY$.
\end{enumerate}

\hfill $\Box$


\vskip 0.7cm

\noindent{\large {\bf Acknowledgments.}} The authors wish to express their thanks to Yuri Sachkov for very fruitful dicussions.


\end{document}